\documentclass[a4paper,11pt]{article}  
\usepackage[utf8]{inputenc} \usepackage{mathtools,amsmath,graphicx,array,amsthm,amssymb}  
\usepackage{float}  
\usepackage{subcaption} 
\usepackage[graphicx]{realboxes}  
\usepackage{hyperref} 
\usepackage{microtype}
\usepackage{bm}

\newtheorem{eg}{Example}[section]  
\let\oldeg\eg  
\let\oldendeg\endeg  
\renewenvironment{eg}{%
    \oldeg  
    \normalfont}{%
    \oldendeg} 

\usepackage{authblk}

\begin{document}

\title{MscaleFNO: Multi-scale Fourier Neural Operator Learning for Oscillatory Function Spaces}
\author[1]{Zhilin You}
\author[1]{Zhenli Xu}
\author[ 2]{Wei Cai\thanks{Corresponding author, email: cai@smu.edu.}}
\affil[1]{School of Mathematical Sciences, MOE-LSC and CMA-Shanghai, Shanghai Jiao Tong University, Shanghai, China}
\affil[2]{Department of Mathematics, Southern Methodist University, Dallas, TX, USA}

\maketitle
\begin{abstract}
    In this paper, a multi-scale Fourier neural operator (MscaleFNO) is proposed to reduce the spectral bias of the FNO in learning the mapping between highly oscillatory functions, with application to the nonlinear mapping between the coefficient of the Helmholtz equation and its solution. The MscaleFNO consists of a series of parallel normal FNOs with scaled input of the function and the spatial variable, and their outputs are shown to be able to capture various high-frequency components of the mapping's image. Numerical methods demonstrate the substantial improvement of the MscaleFNO for the problem of wave scattering in the high-frequency regime over the normal FNO with a similar number of network parameters.
\end{abstract}

\section{Introduction} Operator learning by neural networks has found many applications in establishing cause-effect relations between physical quantities in various physical systems, such as the material properties of a media and its scattering characteristics in wave scattering or the seismic excitation and the response of a building.  Some well-known operator learnings include FNO \cite{li2021fourierneuraloperatorparametric}, DeepOnet 
\cite{Lu_2021_DeepONet,chen1995}, and U-Net \cite{Ronneberger2015UNetCN}. In many situations, the quantities involved may contain many high-frequency contents, which poses a challenge to the deep neural network (DNN) learning due to the well-known spectral bias of the DNN toward low-frequency learning during the training. Various approaches have been proposed to remove the spectral bias of the neural networks such as the multiscale DNN (MscaleDNN) \cite{Ziqi_Liu_2020,Zhang2023ACA} and Phase shift DNN \cite{cai2020}.  Most of the approaches are for the direct learning of the PDE solution. It is our objective in this work to address the spectral bias of operator learning. Specifically, we will extend the MscaleDNN method to the case of the FNO.
There are some existing work on reducing the spectral bias of operator learning by neural networks.  A hierarchical attention neural operator \cite{HANO2024} was proposed to address the spectral bias issue. 
In \cite{oommen2024}, a diffusion model is integrated with various neural operators to improve their spectral representation of turbulent flows.

The rest of the paper is organized as follows. Section 2 presents some preliminaries on FNO and their parameter descriptions. Multi-scale approach for nueral network and the FNO will be introduced in Section 3 as well as the scattering wave mapping from scatter properties will be described. Section 4 contains numerical tests for a simple nonlinear mapping which mimics the Green's function of the Helmholtz equation and the Helmholtz solution mapping. Finally, a conclusion is given in Section 5.

\section{Preliminaries on Fourier Neural Operator}
The Fourier neural operator (FNO) was first proposed by Li {\it et al.} \cite{li2021fourierneuraloperatorparametric} as a data-driven model for learning mappings between infinite-dimensional spaces from observed input-output data. Different from other neural operator architectures such as Deep Green networks \cite{Gin2021DeepGreen,Boull__2022} and DeepONet, the FNO operates by transforming input functions into the frequency domain and learning the operator mapping directly in the spectral space, leveraging the natural structure of the underlying solution space.  

Traditional numerical methods for PDE solving, including finite element and finite difference approaches, require repeated discretization and computation for each new parameter configuration, from different initial conditions, boundary conditions, or physical parameters. This inherent limitation leads to substantial computational cost when addressing  problems with varying parameters. The FNO overcomes this constraint by learning a universal operator that directly maps parameter configurations to their corresponding solutions, thereby enabling efficient inference for novel parameter sets without the need for additional numerical computations.

\subsection{Theory and Architecture}  

Let $ D\subset \mathbb{R}^d$ be a bounded and open set which can be time-dependent. One consider $a \in \mathcal{A}(D; \mathbb{R}^{d_a} )$  and $u\in  \mathcal{U} (D;\mathbb{R}^{d_u})$ as the input and output functions, respectively, where  $ \mathcal{A}(D;\mathbb{R}^{d_a} )$ and $ \mathcal{U} (D;\mathbb{R}^{d_u})$ are spaces of vector function  defined on $D$. 

Let\begin{equation}
    G : \mathcal{A} \rightarrow \mathcal{U}
\end{equation}
is a (nonlinear) mapping such that $G(a) = u$. Suppose we have $N$ observations $\mathbb{S}=\{a_j,u_j\}_{j=1}^{N}$, the goal is to find an approximation $G_\theta$ such that 
\begin{equation}
    G_\theta(a)\approx G(a), \quad \forall a \in  \mathcal{A}(D; \mathbb{R}^{d_a}),
\end{equation}
where $\theta \in\mathbb{R}^{d_p}$ is the parameter of finite dimensions. Thus, our objective can be naturally transformed into solving an optimal problem:
\begin{equation}
\min_{\theta \in \mathbb{R}^{d_p}} \frac{1}{|\mathbb{S}|}\sum_{(a,u)\in\mathbb{S}} L\big(G_\theta(a), u\big),
\label{optimization}
\end{equation}
where the loss function is defined as relative $L_2$ loss
\begin{equation}
L\big(G_\theta(a), u\big) := \frac{\|G_\theta(a)-u\|_{L^2}}{\|u\|_{L^2}}. 
\label{R_L2loss}
\end{equation}
Since both $a$ and $u$ are functions, we need to discretize them for numerical computation. Define $D_n = \{\bm x_1,\cdots, \bm x_n\}$ as the discrete sampling points in the computational domain. In practice, our learning task focuses on constructing the mapping between discrete representations of these functions at a fixed resolution
\begin{equation}
    G: \; \big\{a(\bm x_1),a(\bm x_2),\cdots,a(\bm x_n)\big\}\mapsto \big\{u(\bm x_1),u(\bm x_2),\cdots,u(\bm x_n)\big\}.
\end{equation}
Then, the loss function becomes
\begin{equation}
L\big(G_\theta(a), u\big) := \frac{\sqrt{\sum_{i=1}^n \big( G_\theta(a)(\bm x_i)-u(\bm x_i) \big)^2}}{\sqrt{\sum_{i=1}^n u(\bm x_i)^2}}.
\label{Dis_R_L2loss}
\end{equation}

\noindent\textbf{Neural Operator.} The construction of $G_\theta$ can be achieved through machine learning techniques. Among various approaches, the neural operator framework \cite{BOULLE202483,Kovachki2023NeuralOL} is formulated as an iterative architecture for operator learning, comprising three essential components:
\begin{equation}  
\begin{aligned}  
    &\text{(a)} & v_0(\bm x) & = P(a)(\bm x),\\
    &\text{(b)} &  v_{t}(\bm x) & = \sigma\Big(W_tv_{t-1}(\bm x) + \big(\mathcal{K}(a;\phi_t)v_{t-1}\big)(\bm x)\Big), \\
    &\text{(c)} & u(\bm x) &= Q(v_T)(\bm x).
\end{aligned}  
\label{neural operator}
\end{equation}
The mapping $P : \mathbb{R}^{d_a} \to \mathbb{R}^{d_v}$ represents a linear lifting operator that elevates the input function $a(\bm x)\in \mathbb{R}^{d_a}$ to a higher-dimensional feature space $v(\bm x)\in\mathbb{R}^{d_v}$. In practice, this lifting operator $P$ is typically implemented as a fully connected neural network.  

During the intermediate iteration process, we employ two key operators: a linear local transform $W_t: \mathbb{R}^{d_v} \to \mathbb{R}^{d_v}$, and an integral operator $\mathcal{K}(a;\phi_t): \mathbb{R}^{d_v} \to \mathbb{R}^{d_v}$, which is defined as
\begin{equation}
    \big(\mathcal{K}(a;\phi_t)v_{t-1}\big)(\bm x) := \int_D k(\bm x,\bm y,a(\bm x),a(\bm y);\phi_t)v_{t-1}(\bm y)\,d\bm y.
\end{equation}

The kernel function $k(\bm x,\bm y,a(\bm x),a(\bm y);\phi_t)$ is parameterized by $\phi_t \in \mathbb{R}^{d_\phi}$, with $\sigma$ serving as a nonlinear activation function. The final solution $u(\bm x)$ is obtained through a projection operator $Q: \mathbb{R}^{d_v} \to \mathbb{R}^{d_u}$, typically implemented with a neural network.  
\medskip

\noindent\textbf{Fourier Neural Operator.} The FNO enhances the neural operator framework by imposing translation invariance on the kernel function and eliminating dependence of $a$, resulting in a specialized integral formulation:
\begin{equation}
          \big(\mathcal{K}(a;\phi_t)v_{t-1}\big)(\bm x) 
    = \int_D k_{\phi_t}(\bm x-\bm y)v_{t-1}(\bm y)\,d\bm y.
\end{equation}
By applying the convolution theorem \cite{convolution1998}, we can simplify the integral computation. Let $\mathcal{F}$ denote the Fourier transform operator acting on functions $v: D \to \mathbb{R}^{d_v}$, and $\mathcal{F}^{-1}$ its inverse transform. With $R_t:=\mathcal{F}(k_{\phi_t})$ representing the Fourier transform of the kernel, we obtain
\begin{equation}
   \big(\mathcal{K}(a;\phi_t)v_{t-1}\big)(\bm x) = \mathcal{F}^{-1}\big(R_t\cdot\mathcal{F}(v_{t-1})\big)(\bm x).
\end{equation}

In practical implementations, the output function often contains components that exhibit spatial dependence beyond their direct correlation with the input function. This phenomenon is particularly evident in dynamic systems with source terms, where $u$ represents the solution and $a$ characterizes the initial or boundary conditions or material parameters. In such cases, the solution $u$ inherently incorporates information determined by the spatial distribution of the source terms.  

To effectively handle this spatial dependency, the FNO architecture explicitly incorporates both the spatial coordinates $\bm x$ and the function values $a(\bm x)$ as input variables. Here, the spatial coordinate $x$ plays a crucial role in capturing the position-dependent features of the solution $u(\bm x)$. Mathematically, we can interpret the coordinate input as an identity mapping of $\bm x$:
$$\operatorname{id}_D(\bm x)=\bm x,\quad \forall \bm x\in D \in \mathbb{R}^d. $$
Essentially, FNO learns the mapping
\begin{equation}
    G: \{\operatorname{id}_D\} \times \mathcal{A}(D; \mathbb{R}^{d_a} ) \to \mathcal{U}(D; \mathbb{R}^{d_u} )
\end{equation}
where $u(\bm x) = G\big[\operatorname{id}_D,a\big](\bm x)= G\big[\bm x,a(\bm x)\big](\bm x) $. The comprehensive formulation of the FNO framework can be expressed as follows.
\begin{equation}  
\begin{aligned}  
    &\text{(a)} & v_0(\bm x) & = P\left(\bm x,a(\bm x)\right)(\bm x), \\
    &\text{(b)} &  v_{t}(\bm x) & = \sigma\Big(W_tv_{t-1}(\bm x) + \mathcal{F}^{-1}\big(R_t\cdot\mathcal{F}(v_{t-1})\big)(\bm x)\Big), \\
    &\text{(c)} & u(\bm x) &=  Q\left(v_{T}\right)(\bm x). 
\end{aligned}  
\label{FNO_e}
\end{equation}
$P : \mathbb{R}^{d_a}\times \mathbb{R}^{d} \to \mathbb{R}^{d_v}$ defines an enhanced linear lifting operator that extends the formulation in \eqref{neural operator}, while operators $W_t$ and $Q$ remain consistent with their definitions in \eqref{neural operator}. The spatial coordinate $\bm x$ is incorporated as an optional input parameter, enabling the model to capture position-dependent characteristics in the output function. For notational convenience, we express the equations in \eqref{FNO_e} as

\begin{equation}
    u(\bm x) = \operatorname{FNO}_\theta\big[\bm x,a(\bm x)\big](\bm x).
\end{equation}

\begin{figure}[h]  
\centering  
\includegraphics[height=3.1 cm,width=12cm]{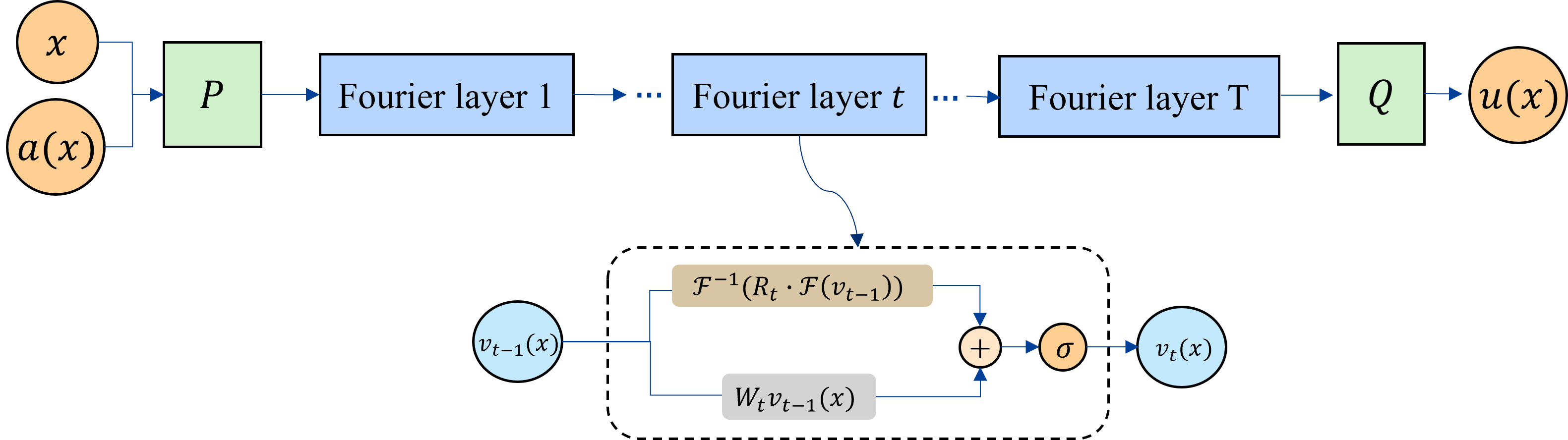}  
\caption{The FNO architecture}  
\label{FNO archi}  
\end{figure}

The comprehensive architecture of the FNO is illustrated in Fig.~\ref{FNO archi}. In the following, we detail the implementation of the truncation mechanism for $R_t$, which constitutes a crucial aspect of the computational framework.

\noindent\textbf{Truncation mechanism of $R_t$.} Considering that  domain $D$ is discretized with $n\in\mathbb{N}$ points such that $v_t\in \mathbb{R}^{d_v\times n}$. When the Fast Fourier Transform (FFT) is applied to $v_t$, it produces $\mathcal{F}(v_t) \in \mathbb{R}^{d_v\times n}$. To achieve computational efficiency, we truncate the Fourier spectrum to retain only the $k_{\mathrm{max}}$ lowest frequency modes, where $k_{\mathrm{max}} < n$.

The truncation operation is implemented through a complex-valued weight tensor $R_t\in \mathbb{C}^{d_v \times d_v \times k_{\mathrm{max}}}$. Mathematically, the truncation process can be expressed as

\begin{equation}  
    \big(R_t\cdot\mathcal{F}(v_{t-1})\big)_{k,l} = \sum_{i=1}^{d_v}(R_t)_{k,l,i}\mathcal{F}(v_{t-1})_{k,i}.
\end{equation}  
where $k=1,\cdots,k_{\mathrm{max}}$ and $i=1,\cdots,d_v$. $R_t$ acts as a masking tensor that preserves the $k_{\mathrm{max}}$ lowest frequency modes and nullify the higher frequencies. This truncation reduces computational complexity while maintaining essential spectral features.

\subsection{FNO Network Parameters}  

To establish a baseline for parameter complexity and facilitate comparison with our proposed model, we give a detailed analysis of the parameter number in the standard FNO model.

The lifting layer $P : \mathbb{R}^{d_a}\times \mathbb{R}^{d} \to \mathbb{R}^{d_v}$ is implemented as a fully connected neural network:  
\begin{equation}  
    P(\bm x,a(\bm x)) = A_{p}^{\bm x} \bm x+ A_{p}^a a(\bm x)+b_p . 
\end{equation}  
where $A_{p}^{\bm x}\in \mathbb{R}^{d\times d_v},\;A_{p}^a\in \mathbb{R}^{d_a\times d_v}$ and $ b_p\in \mathbb{R}^{d_v}$ are parameters and the parameter number for $P$ is:  
\begin{equation}  
 \operatorname{Num}_p =  d\times d_v + d_a\times d_v +d_v = (d+d_a+1)d_v. 
\end{equation}

The Fourier layer consists of two parts. One is the local linear transform  
$W_t(v) = A_{w}v+ b_{w}$  
where $A_{w}\in \mathbb{R}^{d_v\times d_v}$ and $ b_w\in \mathbb{R}^{d_v}$ are parameters. The parameter number for $W_t$ is   
\begin{equation}  
 \operatorname{Num}_w =  d_v\times d_v+d_v = d_v^2+d_v.  
\end{equation}  
The second is The Fourier integral operation with parameter tensor $R_t\in \mathbb{C}^{d_v \times d_v \times k_{\mathrm{max}}}$, which retains $k_{\mathrm{max}}$ modes. The parameter number for this operation is
\begin{equation}  
 \operatorname{Num}_f =  d_v\times d_v\times k_{\mathrm{max}} . 
\end{equation}   

The projection layer $Q: \mathbb{R}^{d_v} \to \mathbb{R}^{d_u}$ is implemented as a two-layer fully connected network, with \hbox{GELU} activation function and $2d_v$ hidden channels:  
\begin{equation}  
    w_T = \hbox{GELU}(A_m v_T+b_m),\quad  
    Q(v_T) = A_{q} w_T+b_q,  
\end{equation}  
where $A_{m}\in \mathbb{R}^{d_v\times 2d_v}$, $A_{q}\in \mathbb{R}^{2d_v\times d_u}$, $ b_m\in \mathbb{R}^{2d_v}$ and $ b_q\in \mathbb{R}^{d_u}$ are parameters. The parameter number is:  
\begin{equation}  
 \operatorname{Num}_q = 2d_v^2+ 2d_v\times d_u +2d_v+d_u.
\end{equation}

The model incorporates a  Multi-Layer Perceptron (MLP) $M: \mathbb{R}^{d_v} \to \mathbb{R}^{d_v}$ following the Fourier integral operation. This modifies \eqref{FNO_e}(b) to:  
\begin{equation}  
    v_{t}(\bm x) = \sigma\bigg(W_tv_{t-1}(\bm x) + M\Big(\mathcal{F}^{-1}\big(R_t\cdot\mathcal{F}(v_{t-1})\big)(\bm x)\Big)\bigg). 
\end{equation}  
The MLP $M$ is implemented as a two-layer fully connected neural network with \hbox{GELU} activation function and hidden dimension $d_v$:  
\begin{equation}  
   M(v) = A_2\big(\text{GELU}(A_1v+b_1)\big)+b_2,\quad \forall\,v \in\mathbb{R}^{d_v},   
\end{equation}  
where $A_1,\,A_2\in\mathbb{R}^{d_v \times d_v}$ are weight matrices and $b_1,\,b_2\in\mathbb{R}^{d_v}$ are bias vectors. The parameter number for the MLP layer is:  
\begin{equation}  
\operatorname{Num}_m = 2d_v^2+2d_v.  
\end{equation}

For the model architecture incorporating MLP after Fourier integral operations, the total number of parameters is:  
\begin{equation}  
\operatorname{Num} = (d+d_a+1)d_v+T[(k_{\mathrm{max}}+3)d_v^2+3d_v]+[2d_v^2+(2d_u+2)d_v+d_u].   
\label{Num_withMLP}  
\end{equation}  

In this architecture, $d_v$ denotes the channel dimension, $k_{\mathrm{max}}$ represents the number of retained Fourier modes, and $T$ indicates the number of Fourier layers. These parameters determine the model size and computational complexity.

\section{Multi-scale FNO}
\subsection{MscaleDNN}
The multi-scale DNN, MscaleDNN,  proposed in \cite{Ziqi_Liu_2020,Zhang2023ACA} uses a simple parallel architectural framework to approximate functions $f:\Omega \to \mathbb{R}$ exhibiting ample spectral content over a general domains $\Omega\subset \mathbb{R}^d$. This architecture is based on the frequency principle \cite{Zhi_Qin_John_Xu_2020} or the spectral bias of DNNs \cite{rahaman2019spectral},  which characterizes the spectral learning dynamics of DNN . The spectral bias of DNNs  refers to the fact that neural networks exhibit preferential learning of low-frequency components over their high-frequency counterparts in the spectral domain. 

To illustrate the main idea of MscaleDNN, let us assume that $f$ is a function with a limited frequency range, which means that its Fourier transform is bounded,   
\begin{equation}  
    \hbox{supp}\hat{f}(\bm k)\subset \{\bm k\in \mathbb{R}^d,\;|\bm k| \leq k_{\mathrm{max}}\} . 
\end{equation}  
First, partitioning the frequency domain of $f(\bm{x})$  
\begin{equation}  
A_i = \{\bm{k}\in \mathbb{R}^d,\; K_{i-1}\leq|\bm{k}| \leq K_i\},\quad i=1,2,\cdots,M,
\end{equation}  
where $0=K_0 < K_1 <\cdots<K_M=k_{\mathrm{max}}$ and $\hbox{supp}\hat{f}(\bm{k})\subset\bigcup_{i=1}^{M} A_i$. Then $f$ can be represented as a sum of functions with non-overlapping frequency information  
\begin{equation}  
f(\bm{x}) = \sum_{i=1}^{M}f_i(\bm{x}),\quad f_i(\bm{x}) = \int_{A_i}\hat{f}(\bm{k})e^{i\bm{k}\cdot \bm{x}}\,d\bm{k}.  
\end{equation}  
Next, one performs a radial scaling on frequency  
\begin{equation}  
  \hat{f}^{(scale)}_i(\bm{k})  = \hat{f}_i(\alpha_i\bm{k}),  
\end{equation}  
which corresponds to its physical-space quantity 
$f_i(\bm{x}) = \alpha_i^n f^{(scale)}_i(\alpha_i\bm{x}).$
If the scaling factor is sufficiently large, the scaled function can be transformed into a low-frequency function with supported Fourier transform  
\begin{equation}  
 \hbox{supp}\hat{f}_i^{(scale)}(\bm{k})\subset \Big\{\bm{k}\in \mathbb{R}^d,\;\frac{K_{i-1}}{\alpha_i}\leq|\bm{k}| \leq \frac{K_i}{\alpha_i}\Big\},\quad i=1,2,\cdots,M.  
\end{equation}  

Due to the Fourier principle, $f^{(scale)}_i(\bm{x})$ can be quickly learned by a DNN $f_{\theta_i}(x)$ parameterized by $\theta_i$,   
\begin{equation}  
    f_i(\bm{x}) =\alpha_i^n f^{(scale)}_i(\alpha_i\bm{x}) \sim \alpha_i^n f_{\theta_i}(\alpha_i \bm{x}).  
\end{equation}  
Since $f(\bm{x})= \sum_{i=1}^{M}f_i\left(\bm{x}\right),$ the framework of the MscaleDNN can be represented by  
\begin{equation}
f(\bm{x}) \sim \sum_{i=1}^{M} \alpha_i^n f_{\theta_i}(\alpha_i\bm{x}).
\label{MscaleDNN_e}
\end{equation}


The MscaleDNN employs a collection of sub-networks with different scaled inputs to approximate the target function $f(\bm{x})$ at different frequency ranges.  A key feature of the method lies in its ability to capture features across multiple scales simultaneously through its parallel network structure. $f_{\theta_i}(\alpha_i\bm{x})$ denotes an individual neural network with parameters $\theta_i$ operating on scaled input $\alpha_i\bm{x}$, where $\alpha_i$  represents the scaling factor. A  property of this structure is that larger values of $\alpha_i$ enable the corresponding sub-network to capture and learn higher-frequency components of the target function $f(\bm{x})$.

\subsection{Multi-scale Fourier Neural Operator (MscaleFNO)}

In complex physical systems, high-frequency components play a vital role in determining system dynamics, yet traditional neural networks often struggle to capture these features accurately. While MscaleDNN has demonstrated remarkable capability in handling multi-scale features for function approximation tasks, the challenges in operator learning are substantially more complex, requiring the capture of high-frequency patterns across both spatial dimensions and function spaces.  

Motivated by these challenges and inspired by MscaleDNN's  approach, we will extend the multiscale concept  to operator learning. Specifically, for the operator mapping:   
$$u(\bm x) = G\big[\operatorname{id}_D,a\big](\bm x)= G\big[\bm x,a(\bm x)\big](\bm x),\quad \bm x\in D.$$  
Our architecture will address two critical aspects: the spatial high-frequency variations with respect to the coordinate position $\bm x$, and the high-frequency oscillations of the operator response to variations in the input function $a(\bm x)$. This dual capability is  significant in physical systems, where the solution dynamics are fundamentally shaped by both local spatial variations and rapid changes in input conditions.  

\begin{figure}[h]  
\centering 
\includegraphics[height=5 cm,width=12 cm]{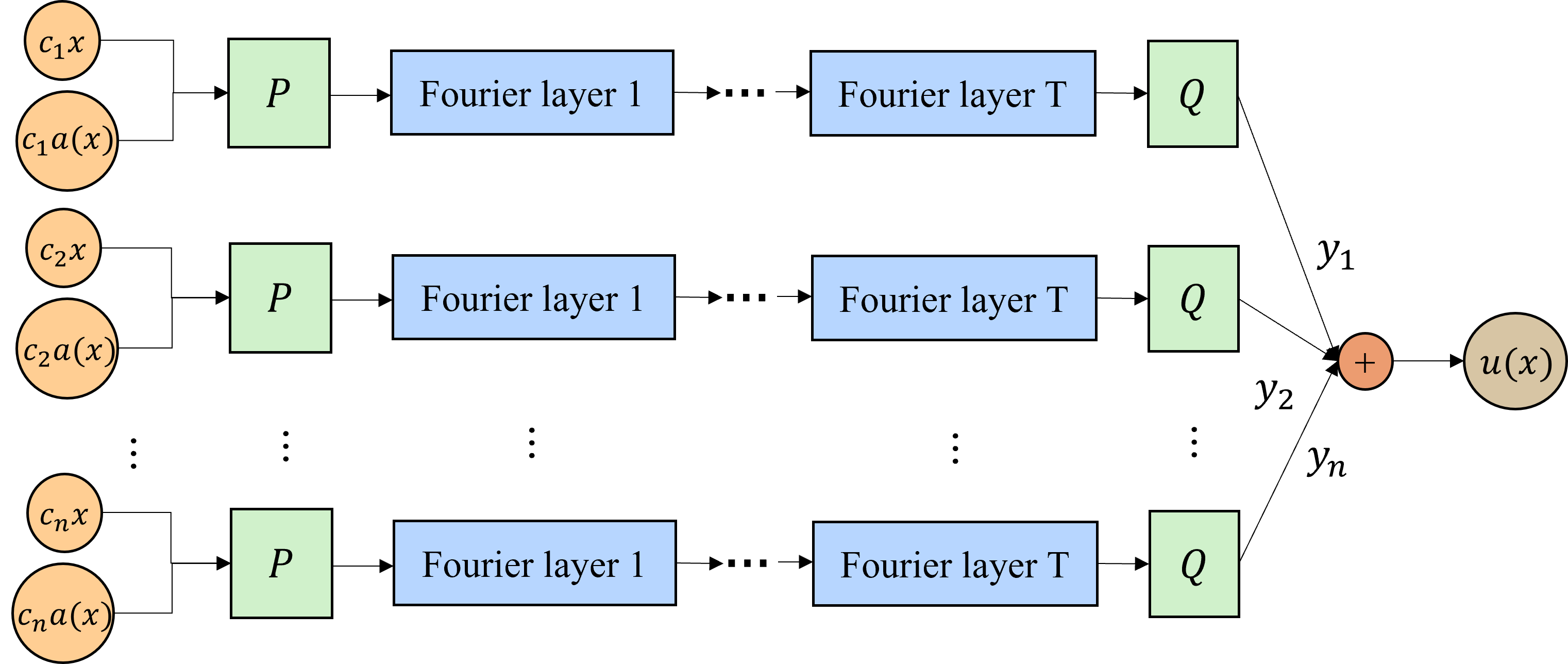}  
\caption{The MscaleFNO architecture}  
\label{MscaleFNO}  
\end{figure}

To achieve this, we propose the following MscaleFNO, as illustrated in Fig.~\ref{MscaleFNO}, a  multi-scale architecture that enhances the Fourier neural operator framework.   
The general mathematical expression of the MscaleFNO model is as follows,  

\begin{equation}  
u\left(\bm x\right) = \sum_{i=1}^{N} \gamma_i \operatorname{FNO}_{\theta_m}\big[c_i \bm x, c_ia(\bm x)\big](\bm x).  
\label{e6}  
\end{equation}  
In this architecture, MscaleFNO employs multiple parallel branches, each implementing a complete FNO structure but processing the input at different scales. These scale parameters simultaneously act on both the spatial variations $\bm x$ and the input field $a(\bm x)$. The final solution is obtained through a weighted summation of outputs from all subnetworks. We denote the standard FNO without this multi-scale architecture as normal FNO for comparison.  

The linear transformations $P$ and $Q$, and the Fourier layers in Fig.~\ref{MscaleFNO} maintain their original definitions from the normal FNO framework. The scaling factors $\{c_i\}_{i=1}^{N}$ for each branch and the weights $\{\gamma_i\}_{i=1}^{N}$ for combining different scales are implemented as trainable parameters. Throughout all Fourier layers, we employ the sine activation function $\sigma(x) = \sin(x).$ 
Given that each sub-network consists of  $d_v$ channels and preserves $k_{\mathrm{max}}$ Fourier modes, the parameter number for each sub-network is $\operatorname{Num}$ defined as Eq. \eqref{Num_withMLP}. 
 The parameter number of MscaleFNO includes both the scaling parameters $\{c_i\}_{i=1}^{N}$ and combination weights $\{\gamma_i\}_{i=1}^{N}$:  
\begin{equation}  
\operatorname{Num}^{(mscale)}=N\times\operatorname{Num}+2N,
\end{equation}  
where $N$ denotes the number of sub-networks in the architecture.  

The multi-scale architecture enables each sub-network to capture frequency components at its corresponding scale, facilitating comprehensive spectral decomposition of the operator. This hierarchical decomposition enhances the precision of operator learning across the frequency spectrum.

\subsection{Mapping between conductivity and solution in Helmholtz equation}
We will consider the following Helmholtz equation for the scattering field of a media with conductivity $a^2(\bm{x})$ under the excitation of incident wave $f(\bm{x})$, 
\begin{equation}
    \Delta u + a^2(\bm{x}) u =f(\bm{x}), \qquad \bm{x} \in \Omega \subset \mathbb{R}^d,
    \label{helm}
\end{equation}
with a boundary condition 
\begin{equation}
    u|_{\partial \Omega} = g(\bm{x}), \quad \bm{x} \in \partial \Omega,
\end{equation}
and the support of $a(\bm{x})$, corresponding to the scatterer, is compact-supported inside $\Omega$.

For a homogeneous boundary condition $g(\bm{x})=0$, the solution of Eq. \eqref{helm} can be written in term of the Green's function $G(\bm{x},\bm{x}^\prime)$,
\begin{equation}
    u(\bm{x})=\int_\Omega G(\bm{x},\bm{x}^\prime) f(\bm{x}^\prime) d \bm{x}^\prime
    \label{ie}
\end{equation}
where the Green's function is defined by
\begin{equation}
    \Delta G(\bm{x},\bm{x}^\prime)  + a^2 G(\bm{x},\bm{x}^\prime) =-\delta(\bm{x}, \bm{x}^\prime),\quad G(\bm{x},\bm{x}^\prime)|_{\bm x\in \partial\Omega}= 0.
    \label{gfe}
\end{equation}

The Green's function $G(\bm{x},\bm{x}^\prime)$ for the homogeneous Dirichlet boundary condition can be written as the sum of the free space Green's function $G_0(\bm{x},\bm{x}^\prime)$ satisfying the Sommerfeld radiation condition \cite{caibook}, plus a smooth function $h(\bm{x}, \bm{x}^\prime)$, i.e.,
\begin{equation}
    G(\bm{x},\bm{x}^\prime) =G_0(\bm{x},\bm{x}^\prime) +h(\bm{x},\bm{x}^\prime)
    \label{gf}
\end{equation}
where
\begin{align}
     G_0(x,x^\prime) & = \frac{i}{2a}e^{i a |x-x^\prime|}, \quad  x, x^\prime \in R, \\
    G_0(\bm{x},\bm{x}^\prime) & = \frac{i}{4}H_0^{(2)}(a|\bm{x}-\bm{x}^\prime|), \quad  \bm{x}, \bm{x}^\prime\in R^2,\\
    G_0(\bm{x},\bm{x}^\prime) & = \frac{e^{i a|\bm{x}-\bm{x}^\prime|}}{4\pi|\bm{x}-\bm{x}^\prime|}, \quad  \bm{x}, \bm{x}^\prime \in R^3,
\end{align} 
and the smooth function $h(\bm{x},\bm{x}^\prime)$ satisfies the following homogeneous Helmholtz equations with non-homogeneous Dirichlet boundary conditions,
\begin{equation}
    \Delta h(\bm{x},\bm{x}^\prime)  + a^2 h(\bm{x},\bm{x}^\prime) =0, \quad h(\bm{x},\bm{x}^\prime)=-G_0(\bm{x},\bm{x}^\prime),\; \bm{x} \in \partial \Omega.
    \label{gfrf}
\end{equation}

\section{Numerical Results}

In the following sections, we present a comparative analysis between the MscaleFNO and the normal FNO with various examples.
To demonstrate the accuracy and efficiency of the MscaleFNO, we make sure that both models have approximately the same number of parameters, and both models use the Adam optimizer with a learning rate of 0.001 for training.
For 1-D function approximation problems where $d_a = d_u = d = 1$, we examine the following two architectures:
\begin{enumerate}
    \item The normal FNO with 48 channels and 500 Fourier modes ($d_v = 48$, $k_{\mathrm{max}} = 500$).
    \item The MscaleFNO with 8 parallel sub-networks, each containing 16 channels and 500 Fourier modes ($d_v = 16$, $k_{\mathrm{max}} = 500$).
\end{enumerate}

\subsection{Learning a nonlinear mapping  between function $a(x)$ and $u(x)$ - $u=\sin m a$.}
In this section, we present numerical results to show the capability of MscaleFNO in capturing high-frequency components of the mapping $u=G[a(x)]$ with respect to the variable $a(x)$. Here, $u(x)$ is defined on the interval $[-1, 1]$.
Considering the dependence of the Helmhotlz solution $u$ on the wave number $a$ in Eq. \eqref{gf}  in 1-D and 3-D cases, we will investigate the following type of nonlinear mapping
\begin{equation}
    u(x)=G[a(x)](x)\equiv \sin (m a(x)) \quad  {\rm or} \quad \cos (m a(x)),
\end{equation}
which will give a mapping of high frequency dependence of $u$  in terms of variable $a$ for large $m$.

The experimental dataset was generated using exact solutions computed on a 1001-point grid. The dataset comprises 2,000 samples,  which were split into 1,000 samples for training, 500 for validation, and 500 for testing. During training, both models utilized a batch size of 20.
Both the MscaleFNO sub-networks and the normal FNO model  employ a single Fourier layer ($T=1$) in their design. Regarding model size,  the MscaleFNO contains {\bf 1,035,544} parameters vs {\bf 1,164,001} parameters in the normal FNO architecture, thus having a smaller number of parameters.
 
\begin{figure}[h]  
\centering  
\subfloat{\includegraphics[width=0.45\linewidth]{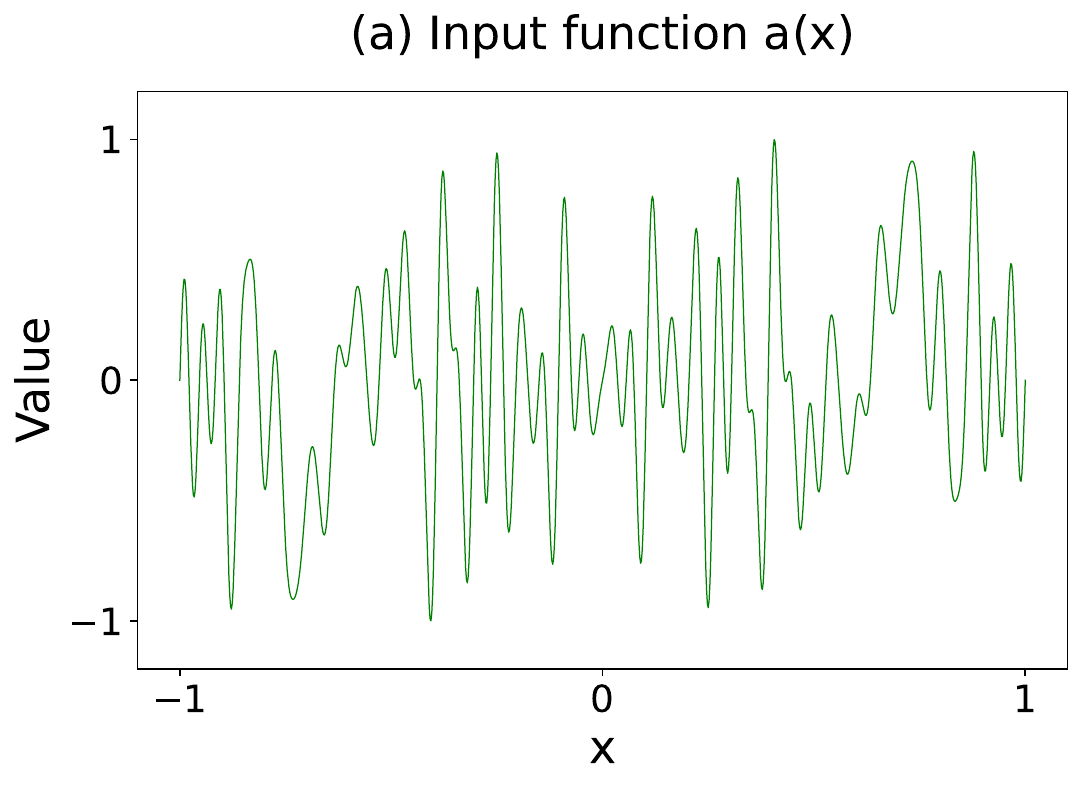}} \hspace{0.04\linewidth} 
\subfloat{\includegraphics[width=0.45\linewidth]{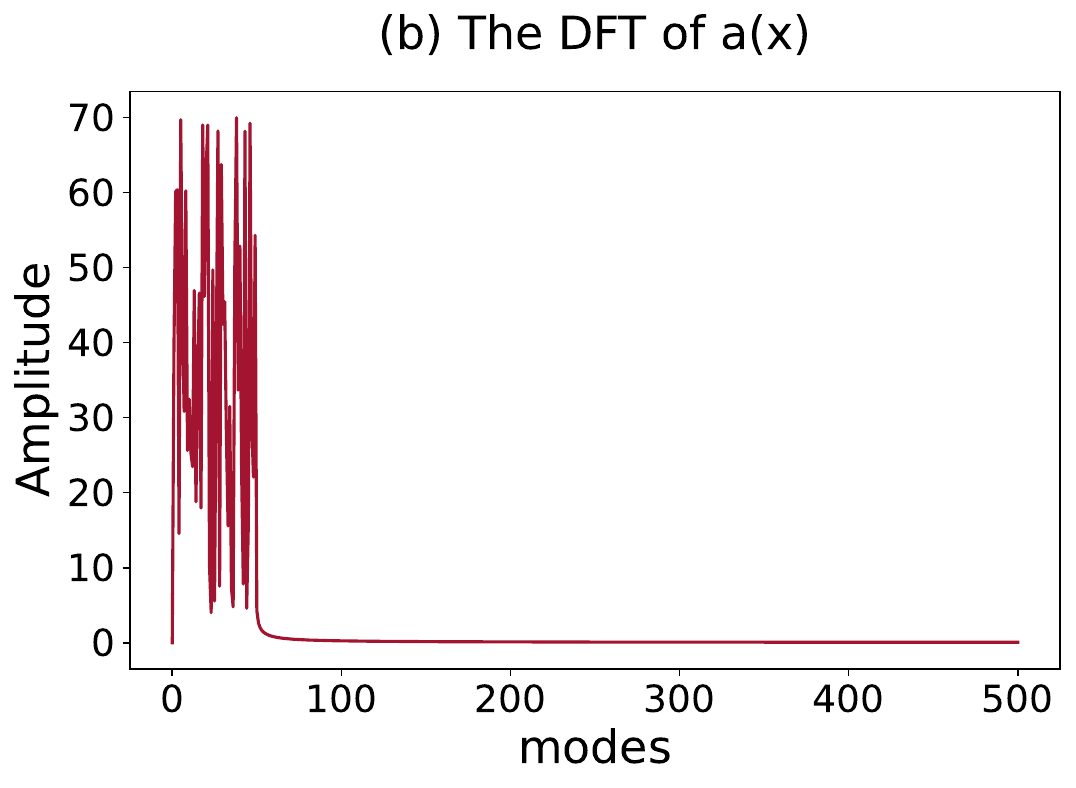}}  
\caption{The profile of the input function $a(x)$ and its DFT} 
\label{a_50}  
\end{figure}

\begin{eg}\textbf{(Single-frequency)}
    We consider the single-frequency function $u(x)$ with respect to $a(x)$, expressed by
\begin{equation}
u(x) = \sin\big(20a(x)\big),\quad x\in [-1,1].
\end{equation}
We aim to learn the operator mapping from $a(x)$ to $u(x)$, where $a(x)$ is generated according to the following form
$$a(x)=\frac{\sum_{n=0}^{50}a_n \sin(n\pi x)}{\max_x\big\{\sum_{n=0}^{50}a_n \sin(n\pi x)\big\}}$$
where $a_n \sim \hbox{rand}(-1,1)$. Fig.~\ref{a_50} shows representative input function $a(x)$ selected from the testing dataset. Denote scales $\bm{c}=\{c_i, i=1, \cdots, 8\}$. We set the initial scaling factors of  MscaleFNO as $\bm{c} = \{1,\;10, \;20, \;40, \;60, \;80, \;100,\;120\}$. 
\end{eg}

\begin{figure}[H]
\centering 
\includegraphics[height=4 cm,width=7 cm]{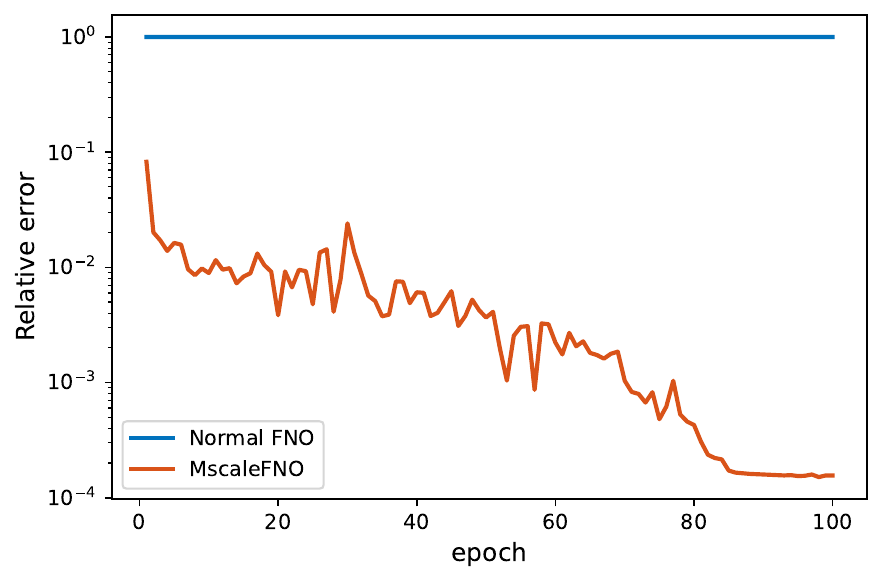}
\caption{Error curves of different models during the training process}
\label{u_sin20a_ec}
\end{figure}

\begin{figure}[h]  
\centering  
\subfloat{\includegraphics[width=0.45\linewidth]{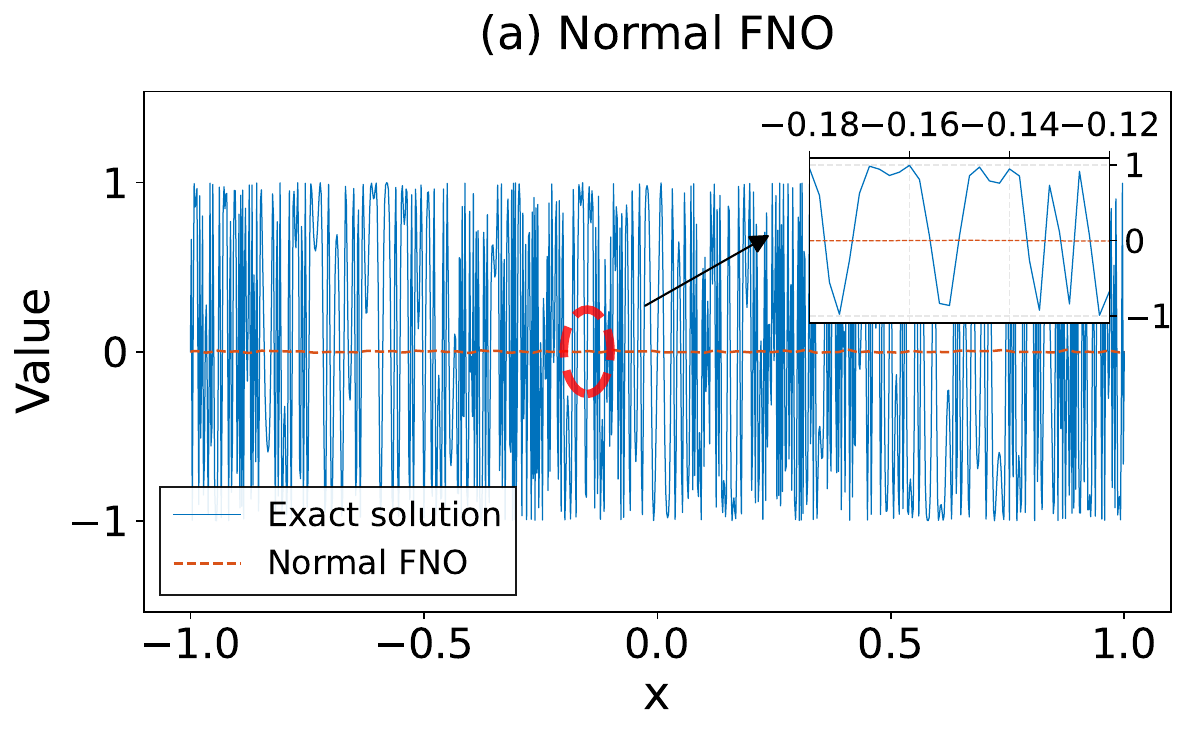}} \hspace{0.04\linewidth} 
\subfloat{\includegraphics[width=0.45\linewidth]{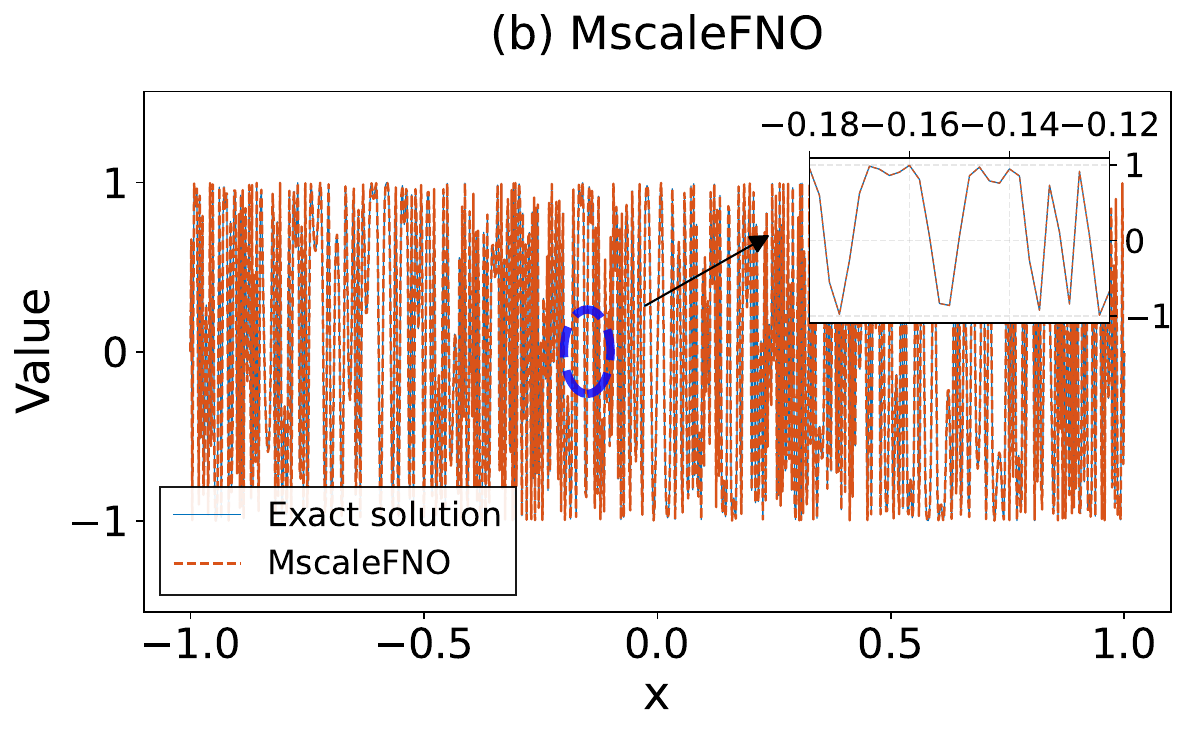}}  
\caption{Predicted solution by the normal FNO (left) and MscaleFNO (right) with zoomed-in inset for $x\in[-0.18,-0.12]$} \label{u_sin20a_spatial}  
\end{figure}

\begin{figure}[h]  
\centering  
\subfloat{\includegraphics[width=0.45\linewidth]{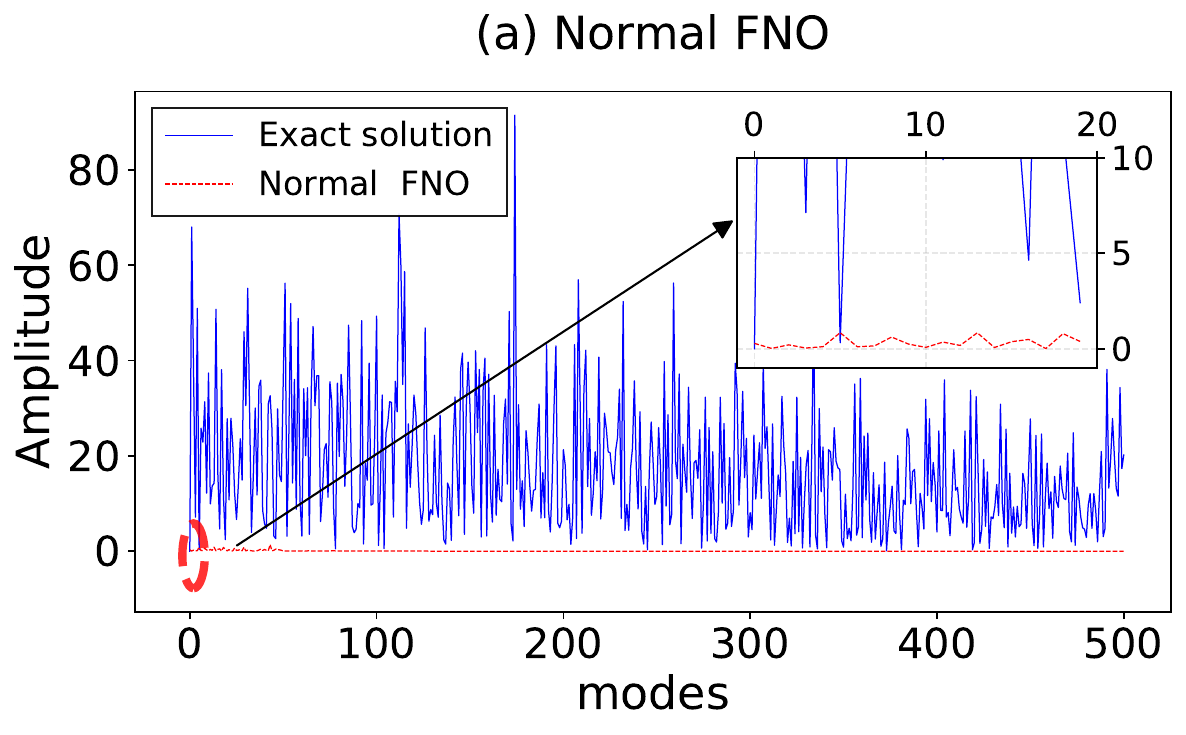}}  
\hspace{0.04\linewidth}
\subfloat{\includegraphics[width=0.45\linewidth]{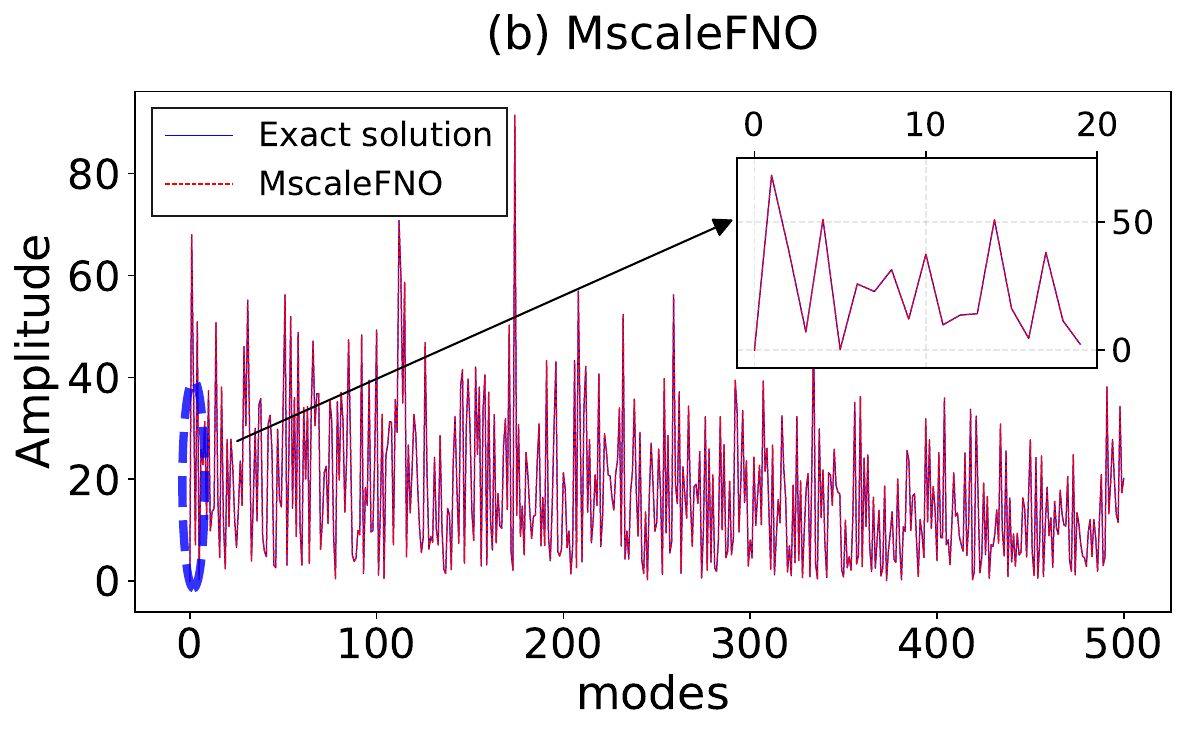}}  
\caption{DFT of predicted solution by normal FNO (left) and MscaleFNO (right) with zoomed-in inset for modes $\in[0,20]$} 
\label{u_sin20a_spec}  
\end{figure}

The relative testing error curves is shown in Fig.~\ref{u_sin20a_ec} , which shows that the test relative error of the MscaleFNO quickly reaches a precision of $O(10^{-4})$ after 100 epochs, while normal FNO stays at $O(1)$. The superior performance of the MscaleFNO is further validated by the predicted solutions shown in Fig.~\ref{u_sin20a_spatial}. The normal FNO fails to capture the high-frequency oscillations in the solution, resulting in a smoothed-out approximation, while the MscaleFNO accurately reproduces the fine wave patterns of the true solution. This observation is quantitatively confirmed by the DFT analysis in Fig.~\ref{u_sin20a_spec}, where the spectrum of the MscaleFNO closely matches the high-frequency components of the true solution, whereas that of the normal FNO shows significant decay in the frequency region. These results collectively demonstrate that  the MscaleFNO exhibits stronger capability in capturing high-frequency components of functions with respect to the input $a(\cdot)$.

\begin{figure}[h]  
\centering  
\subfloat{\includegraphics[width=0.45\linewidth]{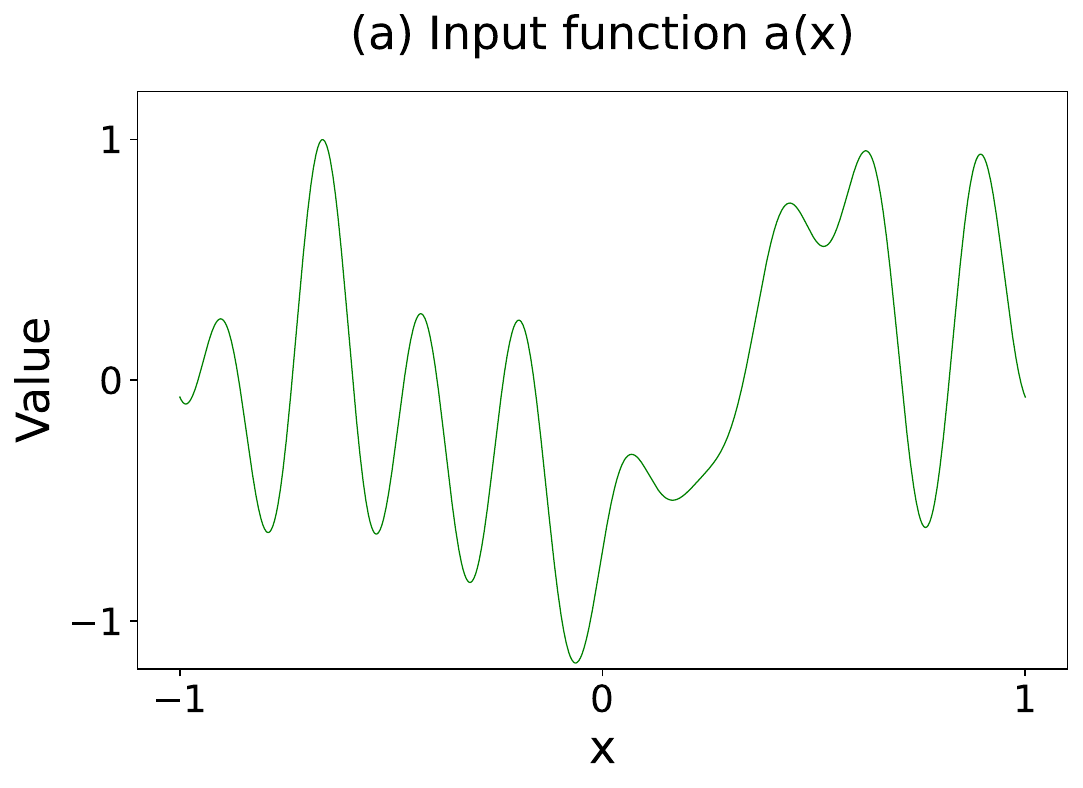}} \hspace{0.04\linewidth} 
\subfloat{\includegraphics[width=0.45\linewidth]{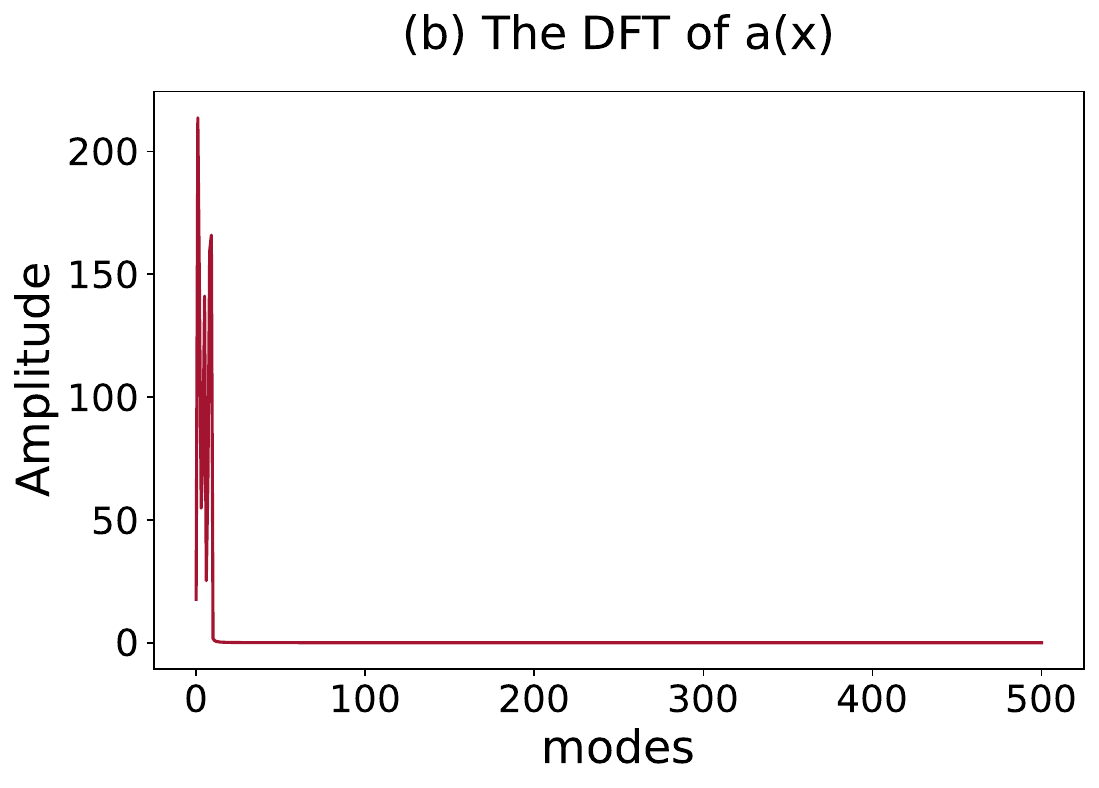}}  
\caption{The profile of the input function $a(x)$ (left) and the DFT of $a(x)$ (right)} 
\label{a_10}  
\end{figure}

\begin{eg}{(\textbf{Multiple-frequency})}
We consider $a(x)$ and $u(x)$ as mixtures of sine and cosine functions, where $u(x)$ exhibits multi-frequency dependence on the input function $a(x)$:
\begin{equation}
u(x) = \sum_{m=1}^{M} \big[A_m\sin \big(ma(x) \big)+B_m\cos \big(ma(x) \big) \big],\quad x\in [-1,1],
\label{u_sum_M}
\end{equation}
where $A_m$ and $B_m\in[-1,1]$ are fixed, and randomly generated numbers in $[-1,1]$.
We are interested in learning the operator mapping from $a(x)$ to $u(x)$ for increasing $M=10,\;20,\;40,\;80,\;100$ and $200$.
The $a(x)$ is generated according to 
$$a(x)=\frac{\sum_{n=0}^{10}[a_n \sin(n\pi x)+b_n \cos(n\pi x)]}{max_x\big\{\sum_{n=0}^{10} [a_n \sin(n\pi x)+b_n \cos(n\pi x)]\big\}}$$ 
where $a_n,b_n \sim \hbox{rand}(-1,1)$. A representative input function $a(x)$ from the testing dataset is visualized in Fig.~\ref{a_10}.
\end{eg}

\begin{figure}[h]  
\centering  
\subfloat{\includegraphics[width=0.3\linewidth]{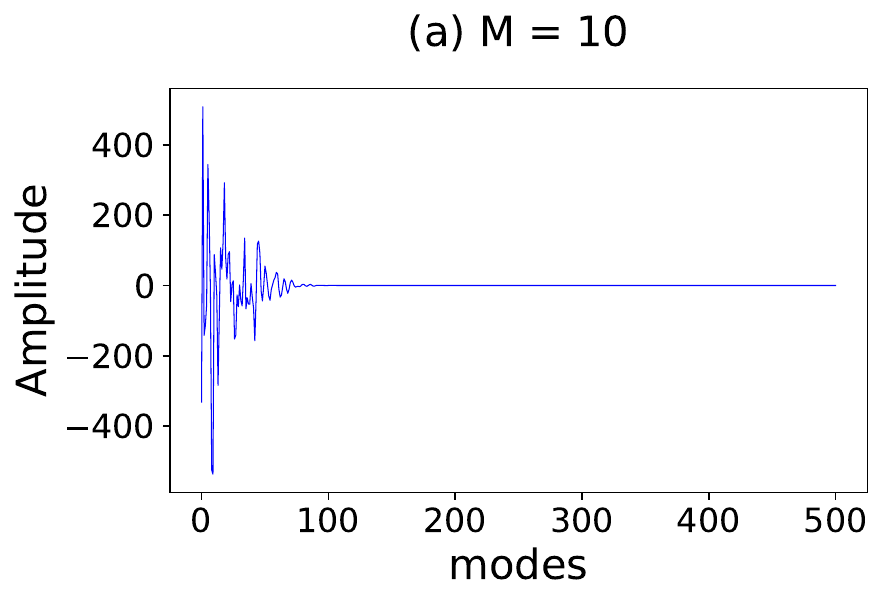}}  
\subfloat{\includegraphics[width=0.3\linewidth]{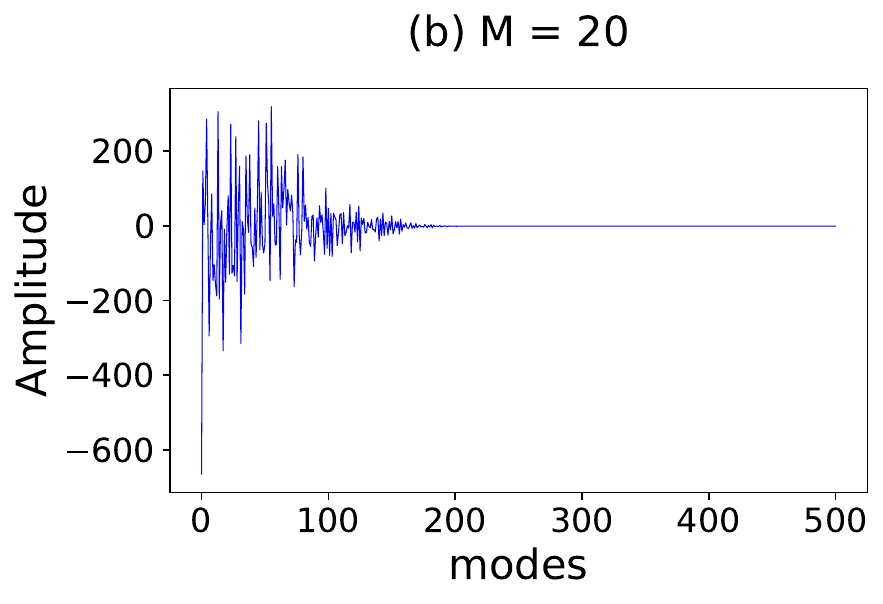}}  
\subfloat{\includegraphics[width=0.3\linewidth]{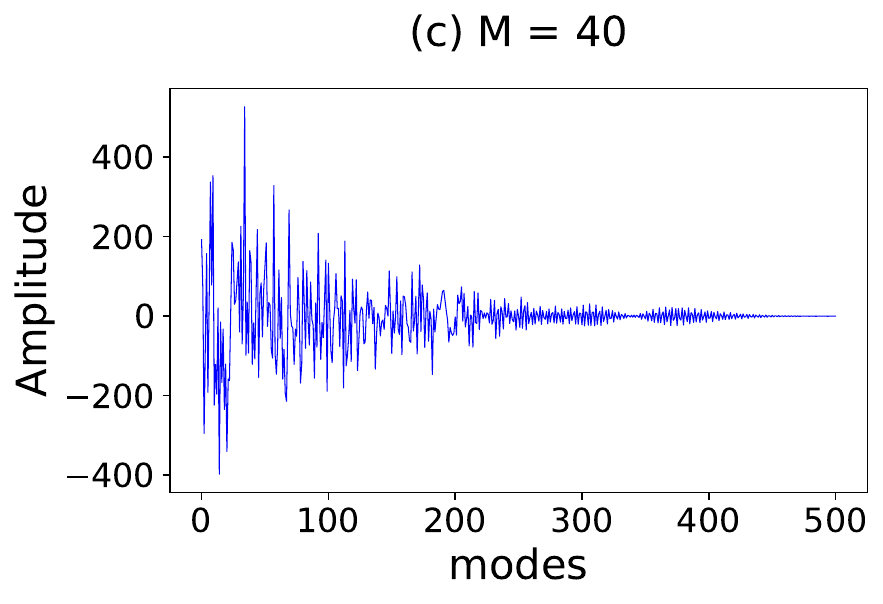}} \\
\subfloat{\includegraphics[width=0.3\linewidth]{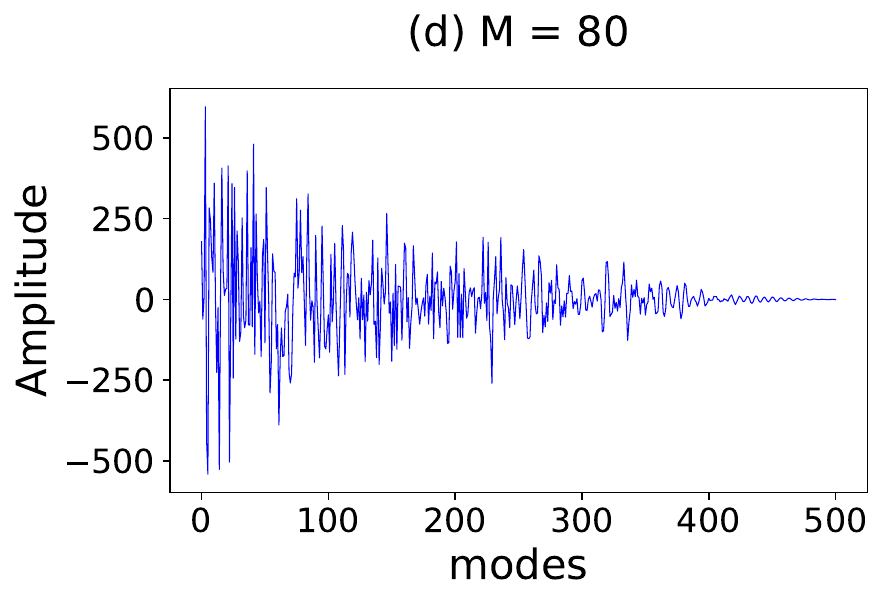}}  
\subfloat{\includegraphics[width=0.3\linewidth]{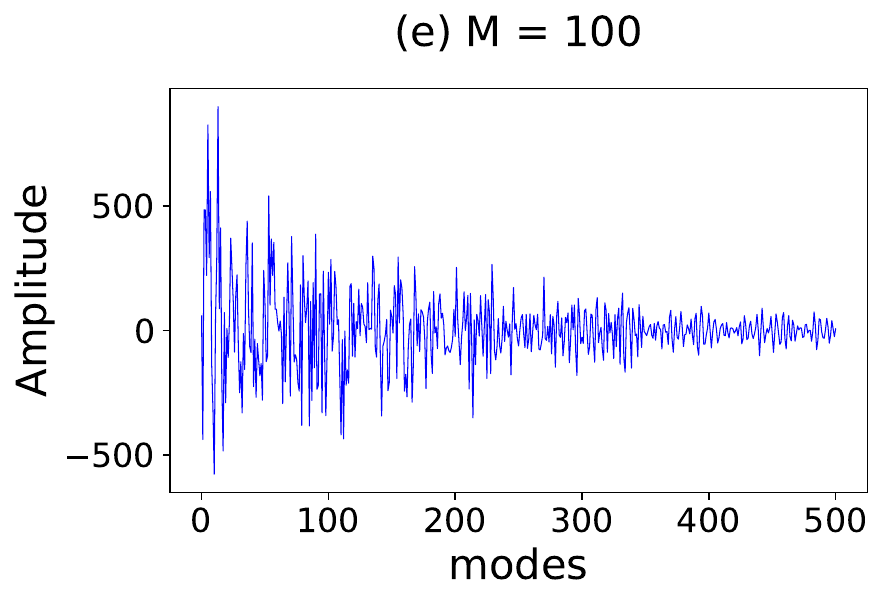}} 
\subfloat{\includegraphics[width=0.3\linewidth]{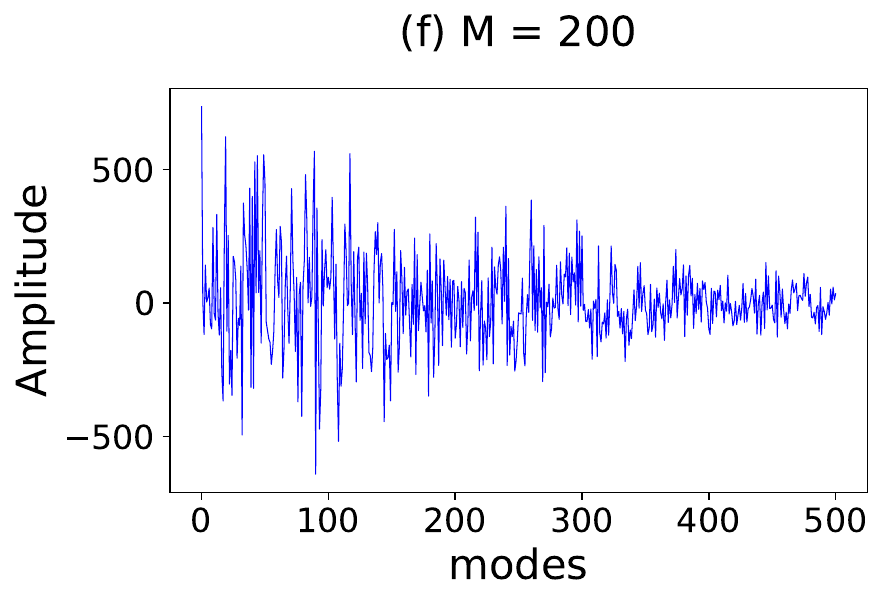}}
\caption{DFT of representative exact solution $u(x)$ for different $M$} 
\label{u_Msum_freq}  
\end{figure}

We visualize its corresponding solution in spectral domains(Fig.~\ref{u_Msum_freq}). As $M$ (the number of frequency terms) increases from 10 to 200, we observe high complexity in the solution's spectral content, particularly in the high-frequency regime, as shown in  their DFT analysis in Fig.~\ref{u_Msum_freq}. In the frequency domain, the spectrum evolves from a concentrated distribution near low frequencies to an increasingly broader distribution extending into higher frequencies as $M$ grows. This spectral distribution manifests itself in the spatial domain as an increasingly oscillatory behavior.

\begin{figure}[h]  
\centering  
\subfloat{\includegraphics[width=0.3\linewidth]{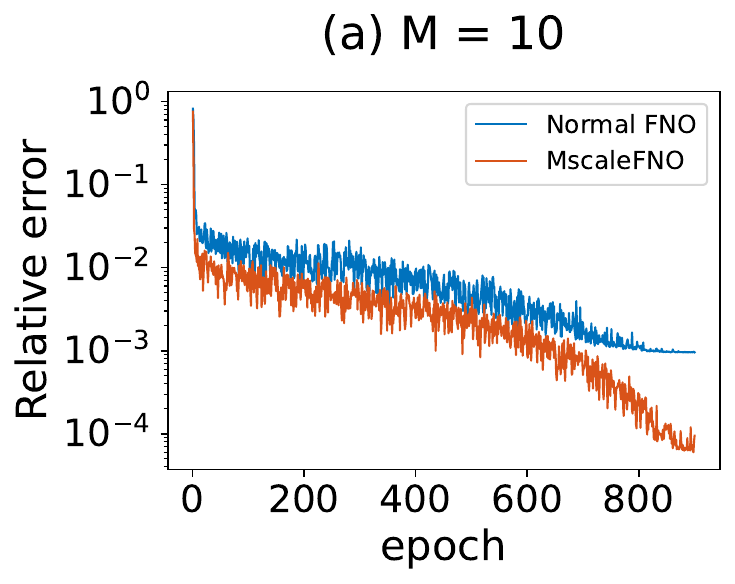}}  
\subfloat{\includegraphics[width=0.3\linewidth]{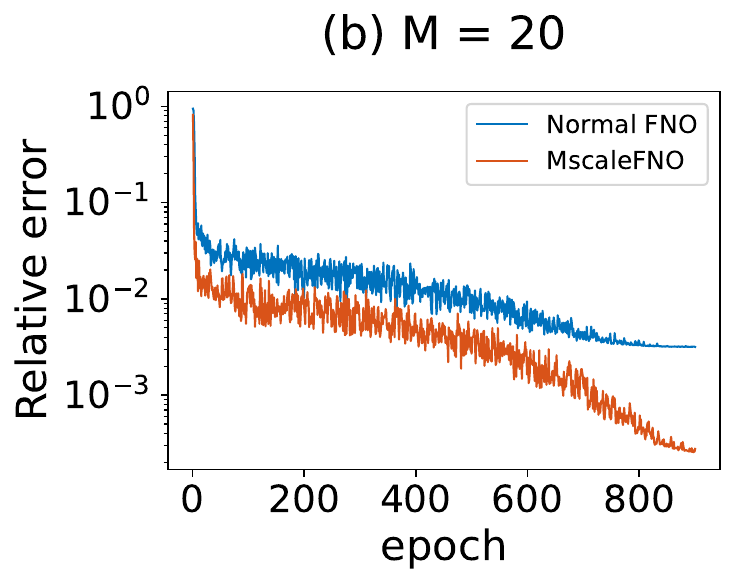}}  
\subfloat{\includegraphics[width=0.3\linewidth]{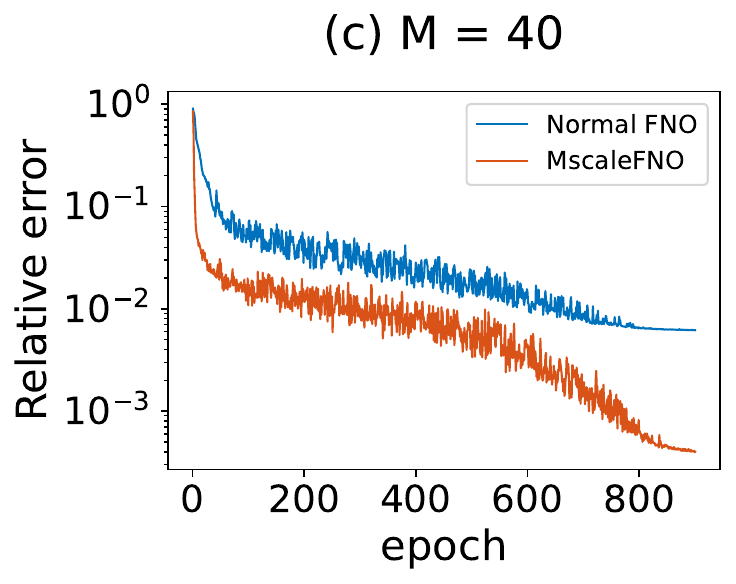}} \\
\subfloat{\includegraphics[width=0.3\linewidth]{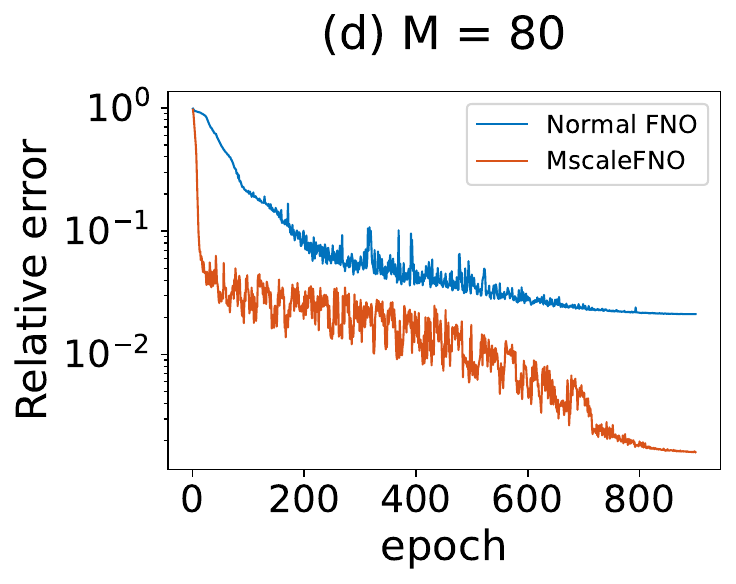}}  
\subfloat{\includegraphics[width=0.3\linewidth]{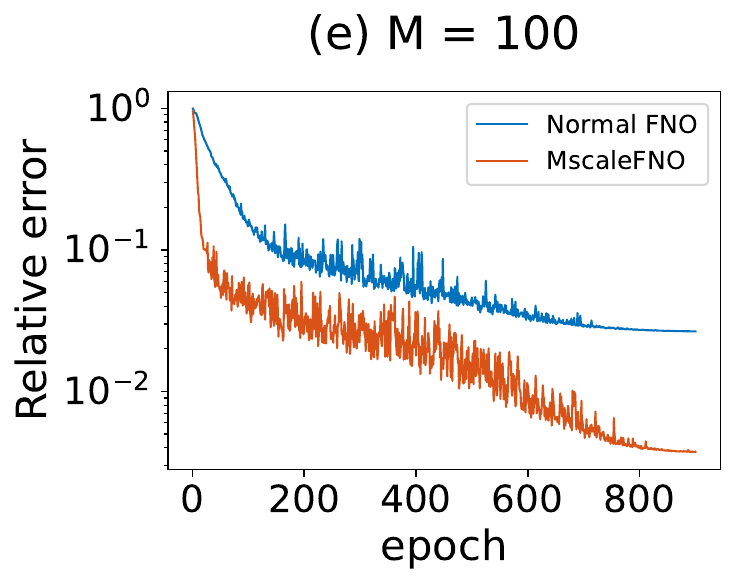}} 
\subfloat{\includegraphics[width=0.3\linewidth]{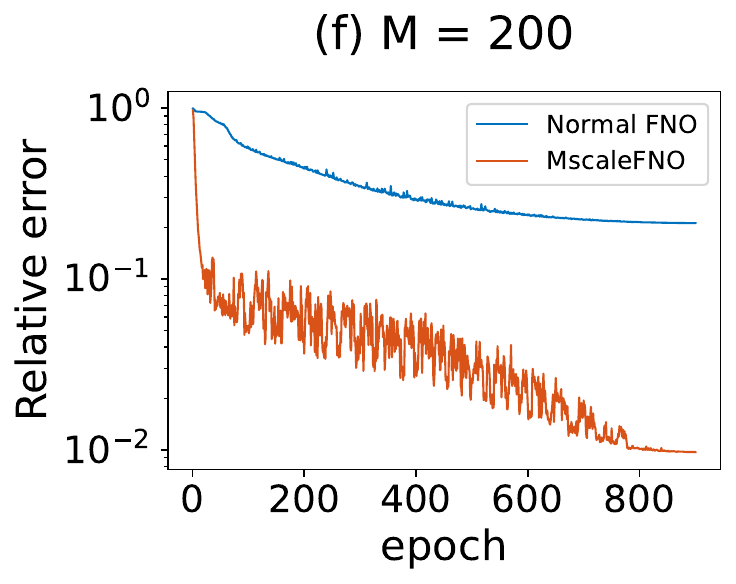}}
\caption{Error curves of different models during the training process under different values of $M$ (Epoch=900)} 
\label{u_Msum_ec}  
\end{figure}

Fig.~\ref{u_Msum_ec} presents the relative testing error curves of the MscaleFNO (blue) and normal FNO (orange) under different frequency modes $M$ (ranging from 10 to 200) over 900 epochs.  For the MscaleFNO, the initial setting of scales in the cases of $M=1,\;20,\;40,\;80,\;100$ takes $\bm{c} = \{1, \;10, \;20, \;40, \;60, \;80, \; 100, \; 120\}$. In the case of $M=200$, it takes $\{1,\; 40,\; 80, \;100,\; 120,\; 140,\; 180, \; 200\}$. We can see that for the normal FNO (blue), as $M$ increases, which corresponds to a scatterer size of more wavelengths, and thus to a high-frequency scattering regime, the relative testing error shows a significant increase. At $ M=200$, the normal FNO has a relative error of $0.2$. In contrast, the MscaleFNO (orange) not only consistently outperforms the normal FNO but also maintains an accuracy of relative error at the level of $10^{-2}$. 

\begin{figure}[h]  
\centering  
\subfloat{\includegraphics[width=0.45\linewidth]{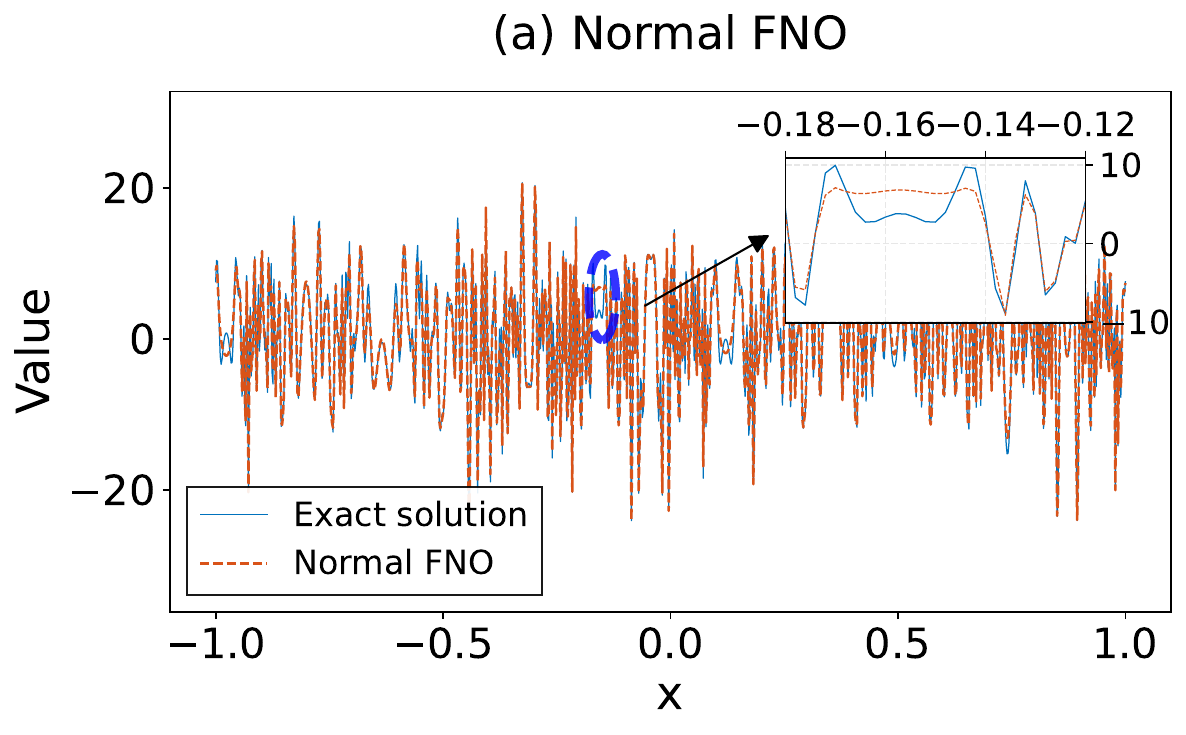}}  
\subfloat{\includegraphics[width=0.45\linewidth]{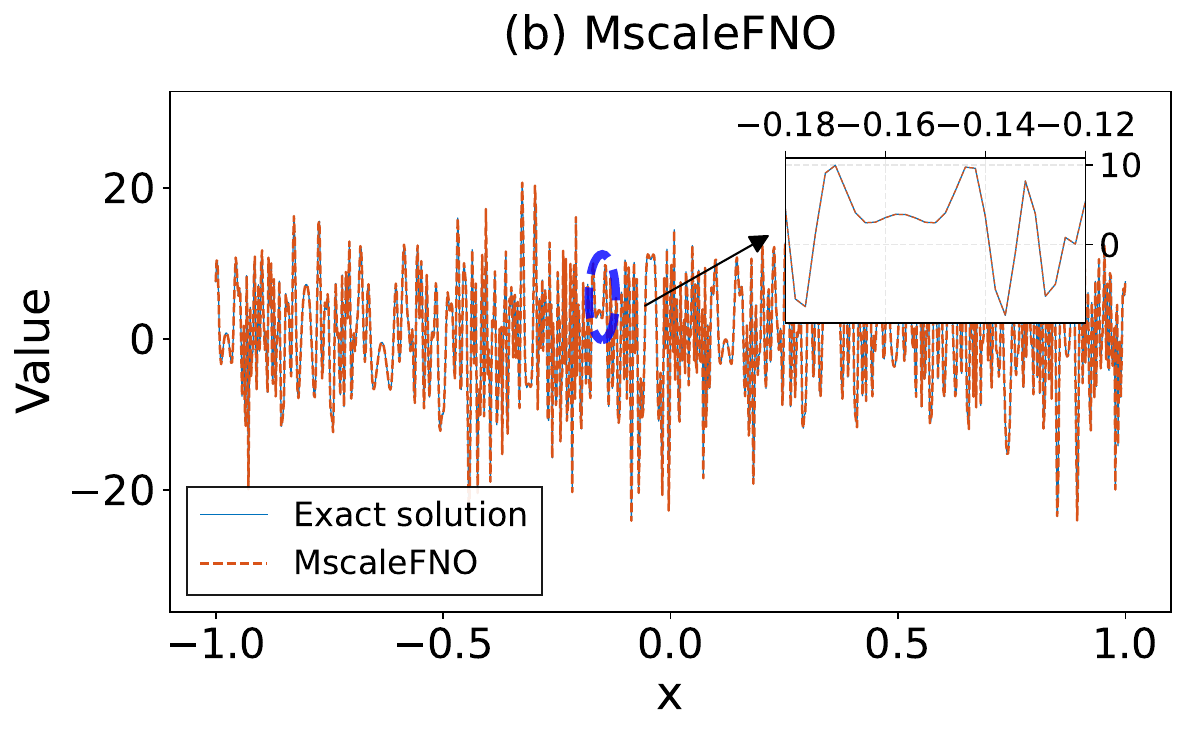}}  
\caption{$M = 200$: Predicted solution by normal FNO (left) and MscaleFNO (right) with zoomed-in inset for $x\in[-0.18,-0.12]$} 
\label{u_sinMa_M200_spatial}  
\end{figure}

\begin{figure}[h]  
\centering  
\subfloat{\includegraphics[width=0.33\linewidth]{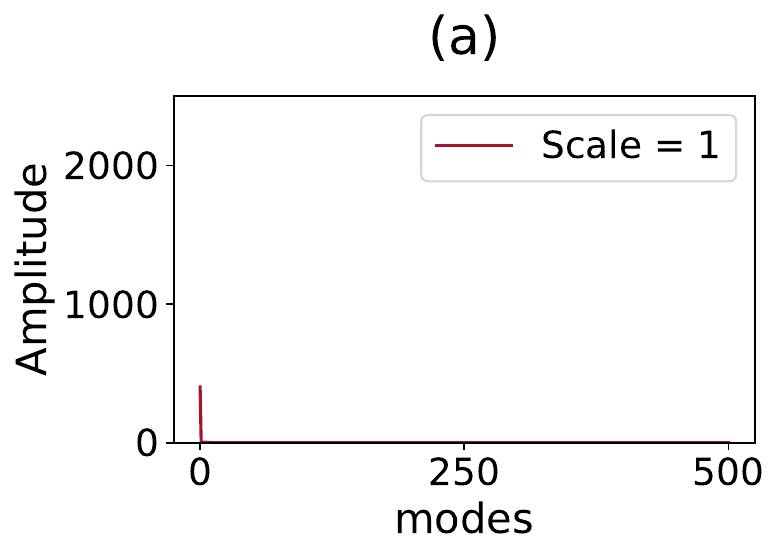}}  
\subfloat{\includegraphics[width=0.33\linewidth]{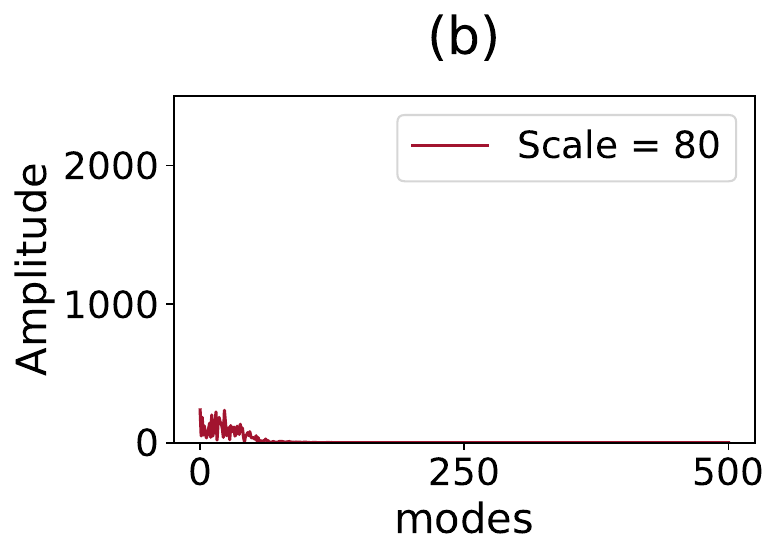}}  
\subfloat{\includegraphics[width=0.33\linewidth]{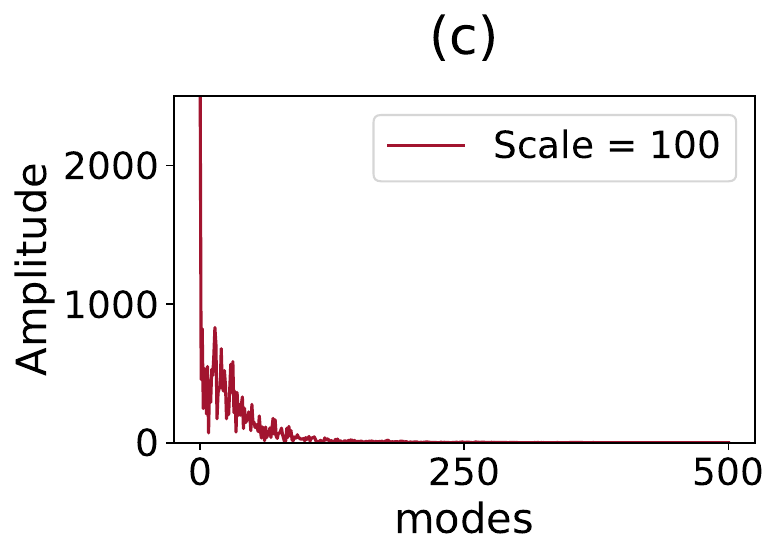}}  \\
\subfloat{\includegraphics[width=0.33\linewidth]{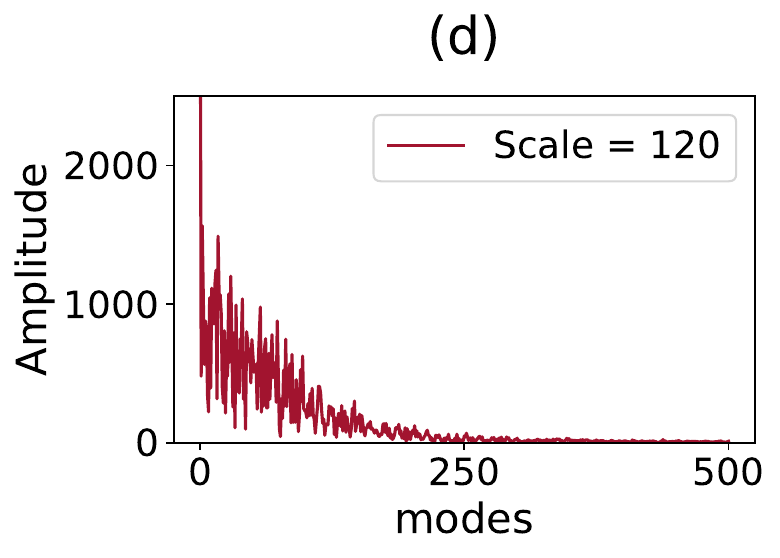}}  
\subfloat{\includegraphics[width=0.33\linewidth]{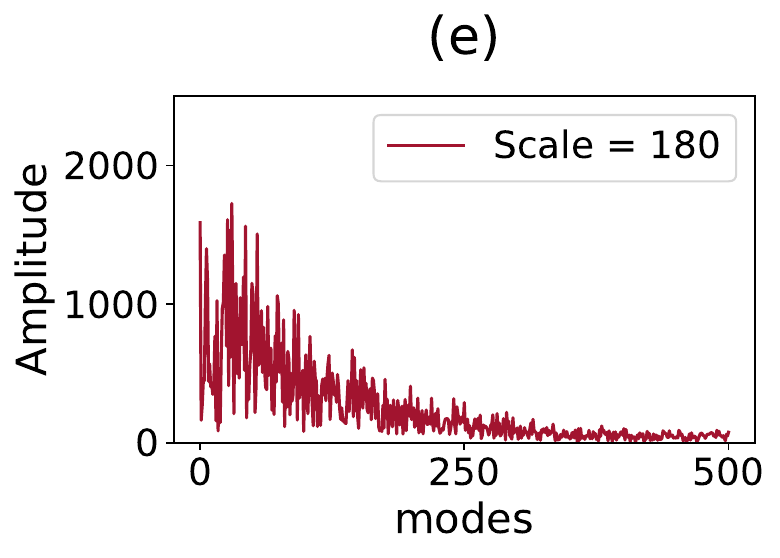}}  
\subfloat{\includegraphics[width=0.33\linewidth]{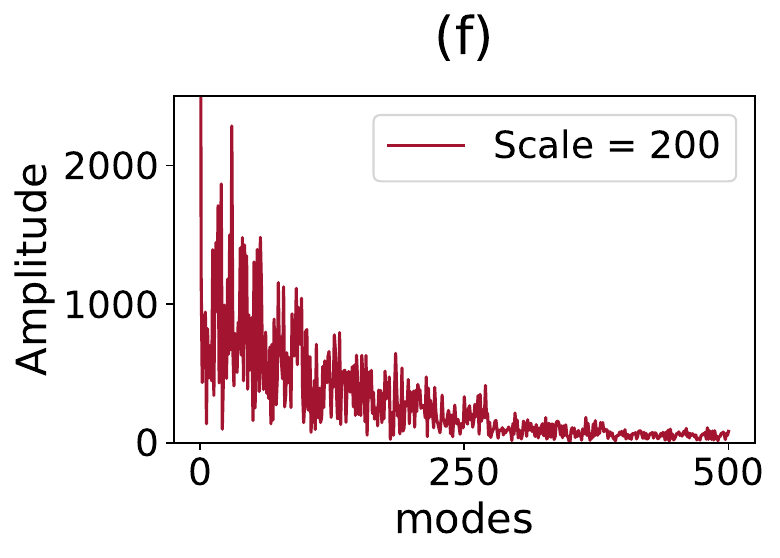}}  
\caption{$M=200$: Spectral contributions of MscaleFNO subnetworks corresponding to differet initial scales} 
\label{u_Msum_M200_outofQ}  
\end{figure}
Fig.~\ref{u_sinMa_M200_spatial} illustrates the predicted solution by the normal FNO and MscaleFNO for the case $M=200$ where the solution displays highly oscillatory patterns. The inserts clearly show the superior performance of the MscaleFNO in capturing the fine details of the solutions.
Fig.~\ref{u_Msum_M200_outofQ} presents the discrete Fourier transform of the contributions $\{\gamma_i\operatorname{FNO}_{\theta_i}\big[c_ix,c_ia(x)\big]\}_{i=1}^8$ from different subnetworks in MscaleFNO at $M=200$, revealing a systematic pattern in the frequency decomposition across different scales. Specifically, the subnetworks with increasing scale parameters exhibit enhanced capability in capturing high-frequency components, as evidenced by their progressively broader frequency spectrum distributions. This hierarchical analysis demonstrates that subnetworks with different $c_i$ values exhibit complementary frequency responses. The subnetworks with smaller $c_i$ capture low-frequency patterns, while those with larger $c_i$ extract high-frequency details, collectively achieving a complete frequency spectrum representation.

\subsection{Mapping between scatterer properties and solution for Helmholtz Equation}
We consider the 1-D Helmholtz equation of (\ref{helm}) with 
$a^2(x)=  \lambda^2 + c\omega(x)$
for a scatterer  with Dirichlet boundary condition, namely,
\begin{equation}
    \begin{cases}  
    u'' + (\lambda^2 + c \omega(x))u = f(x), & x \in [-L,L] \\
    u(-L) = u(L) = 0,  
\end{cases}
\label{e10}
\end{equation}
where $\omega(x)$ is the variable wave number perturbation caused by the scatterer, $u(x)$ is the scattering field and $f(x)$ is the forcing term resulting from an external wave source. We are interested in learning the operator mapping from the variable wave number $\omega(x)$  to the solution $u(x)$. We choose $\lambda = 2$, $c = 0.9\lambda ^2 = 3.6$ and 
\begin{equation}
        f(x) = \sum_{k=0}^{10} (\lambda^2-\mu_k^2) \sin(\mu_k x), \; \mu_k = 300 + 35k.
        \label{e11}
\end{equation}
The perturbation $\omega(x)$ is generated according to 
\begin{equation}
        \omega(x)=\frac{\sum_{n=0}^{500}[a_n \sin(n\pi x)+b_n \cos(n\pi x)]}{\max_x\big\{\sum_{n=0}^{500}[a_n \sin(n\pi x)+b_n \cos(n\pi x)]\big\}}
        \label{e12}
\end{equation}
where $a_n,b_n \sim \hbox{rand}(-1,1)$.

For each problem, we generated a dataset consisting of 1,000 samples, partitioned into 800 for training, 100 for validation, and 100 for testing. Both the normal FNO and MscaleFNO models were trained using a batch size of 20, and  employ four Fourier layers ($T=4$) in their design. With parameter numbers of {\bf 4,127,128} and {\bf 4,641,169}, respectively, and the MscaleFNO has a smaller parameter count.

\begin{figure}[h]  
\centering  
\subfloat{\includegraphics[width=0.45\linewidth]{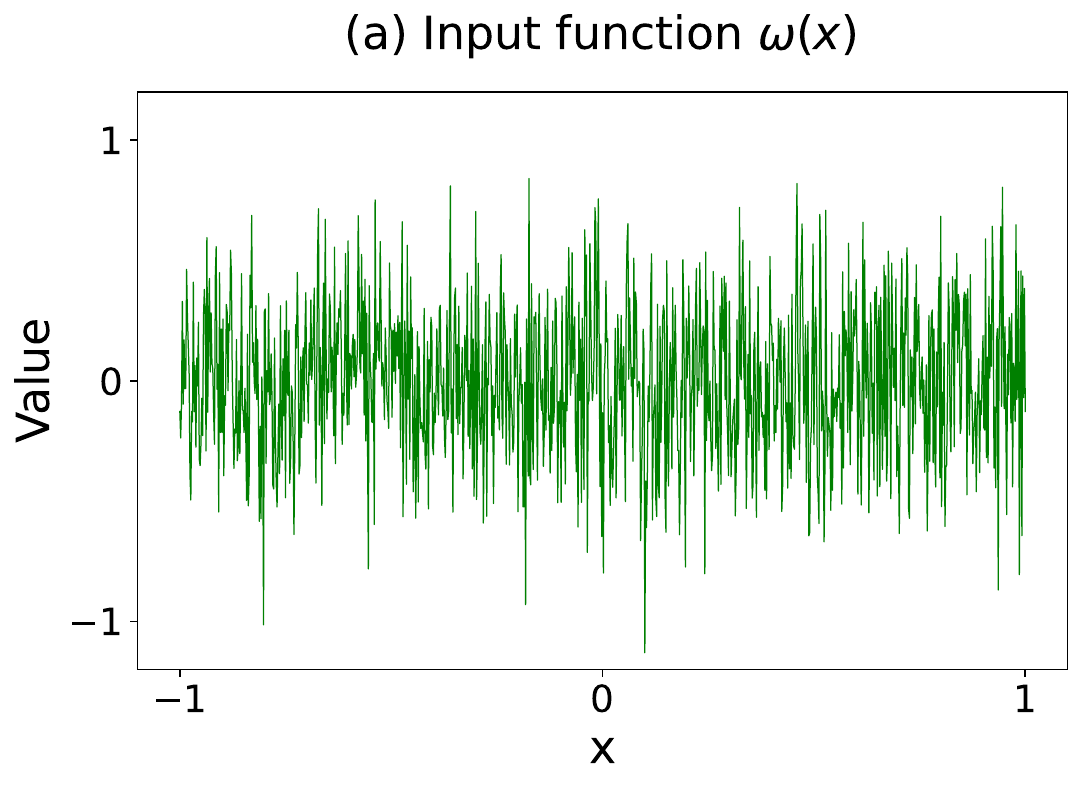}} \hspace{0.04\linewidth}
\subfloat{\includegraphics[width=0.45\linewidth]{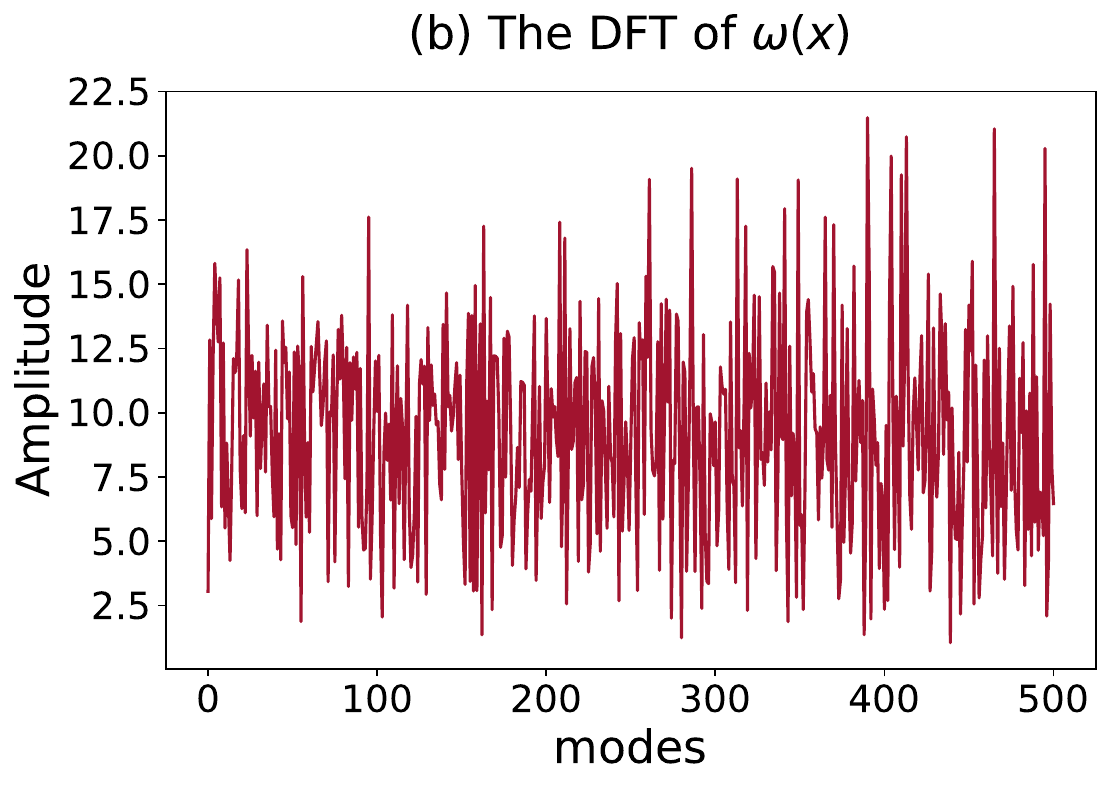}}  
\caption{The profile of the input function $\omega(x)$ (left) and the DFT of $\omega(x)$ (right)} 
\label{helm_w500_w}  
\end{figure}

\begin{figure}[h]
\centering 
\subfloat{\includegraphics[width=0.45\linewidth]{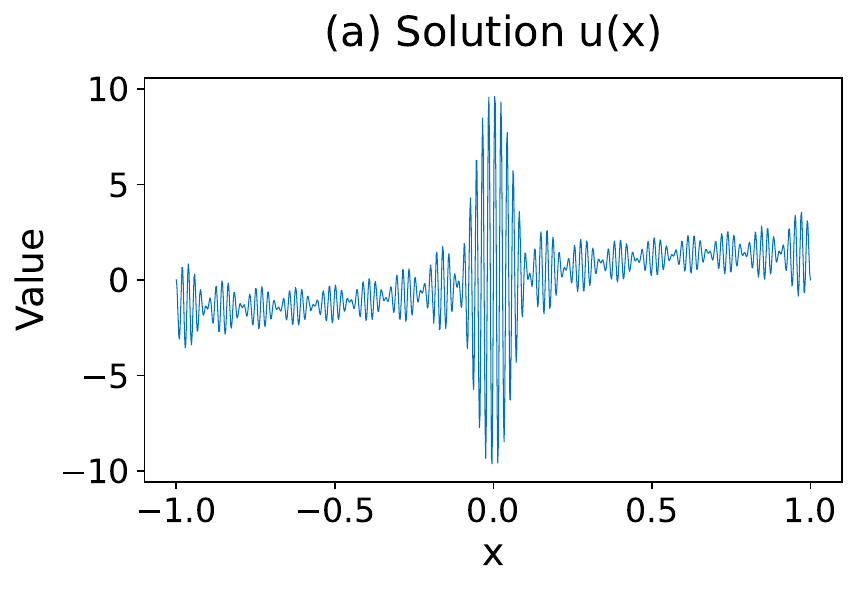}}  
\hspace{0.05\linewidth}
\subfloat{\includegraphics[width=0.45\linewidth]{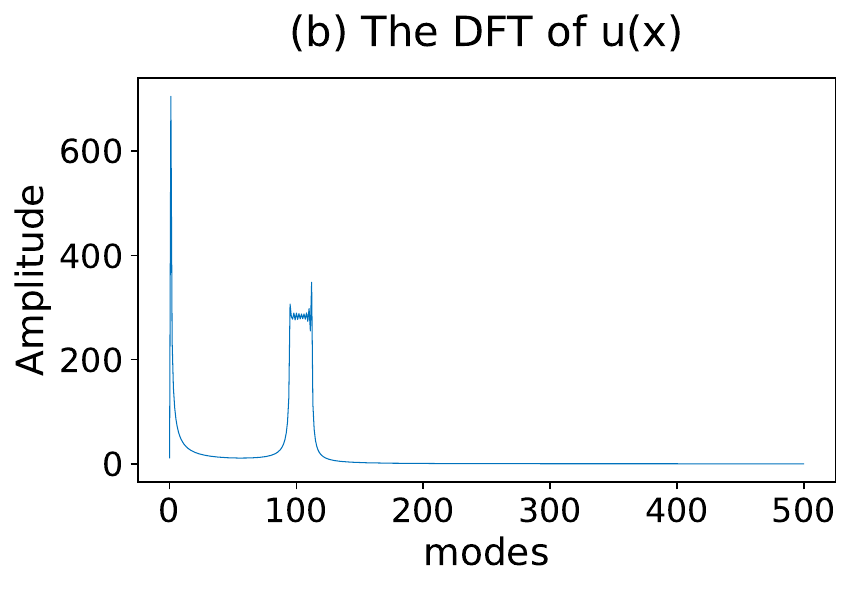}} 
\caption{The profile of the exact solution $u(x)$ from $\omega(x)$ in Fig. \ref{helm_w500_w} (left) and the DFT of $u(x)$ (right)} 
\label{helm_w500_true}
\end{figure}

\begin{eg}
 We first consider the Helmholtz equation Eq. $\eqref{e10}$ with $L=1$. The equation is solved numerically using the finite difference method with a high-resolution grid of 8001 points, achieving a numerical  accuracy of relative error at  $O(10^{-4})$. The high-resolution solution is subsequently downsampled to a coarser grid of 1001 points for computational efficiency. Fig.~\ref{helm_w500_w} illustrates the characteristic profiles of input function $\omega(x)$ and the corresponding solution $u(x)$ is shown in Fig.~\ref{helm_w500_true}.
\end{eg}

\begin{figure}[h]
\centering 
\includegraphics[height=4cm,width=7 cm]{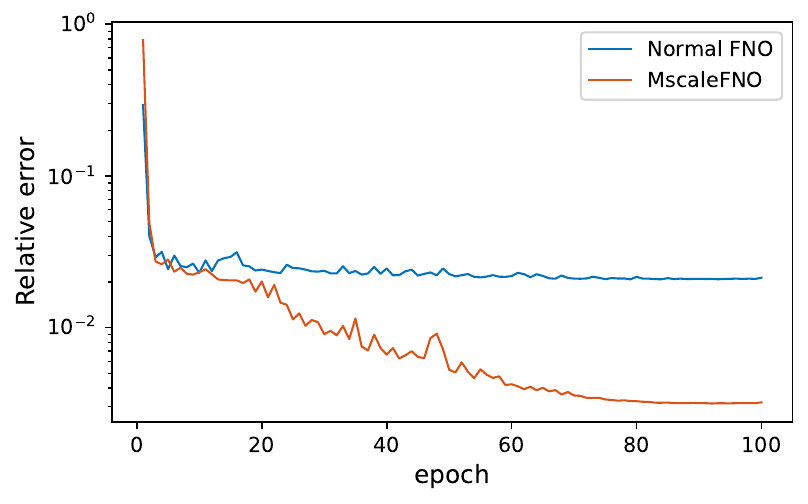}
\caption{Error curves of different models during the training process}
\label{helm_500_o_ec}
\end{figure}

The initial scales of MscaleFNO are set to $ \{1,\; 4,\; 8,\;10,\; 12,\; 14,\; 18,\; 20\}$. Fig.~\ref{helm_500_o_ec} shows the relative testing errors of normal FNO and MscaleFNO over 100 training epochs. Both architectures demonstrate rapid convergence in the initial 10 epochs, followed by different convergence behaviors. The normal FNO converges to a relative error at $O(10^{-2})$ with limited further improvement, while the MscaleFNO continues to reduce errors throughout training, reaching a relative error at $O(10^{-3})$ and achieving an order-of-magnitude improvement in accuracy.

\begin{eg}
We investigate the Helmholtz equation across varying spatial domains characterized by lengths $L \in \{2,4,8,10\}$, corresponding to progressively increasing scatterer size or higher frequency scattering problems. The perturbation $\omega(x)$ is constructed according to
\begin{equation}
        \omega(x)=\frac{\sum_{n=0}^{50}[a_n \sin(n\pi x)+b_n \cos(n\pi x)]}{\max_x\big\{ \sum_{n=0}^{50}[a_n \sin(n\pi x)+b_n \cos(n\pi x)]\big\}}
\end{equation}
where $a_n,b_n \sim \hbox{rand}(-1,1)$.
For each $L$, Eq.~$\eqref{e10}$ is solved numerically using the finite difference method with sufficiently refined grids, ensuring a numerical accuracy of $O(10^{-4})$ relative error for all cases. To maintain consistent spatial resolution with the case ($L=1$), these high-resolution solutions are subsequently downsampled to achieve a uniform mesh size, resulting in different numbers of grid points proportional to the respective domain lengths.
\end{eg}

\begin{figure}[h]  
\centering  
\subfloat{\includegraphics[width=0.35\linewidth]{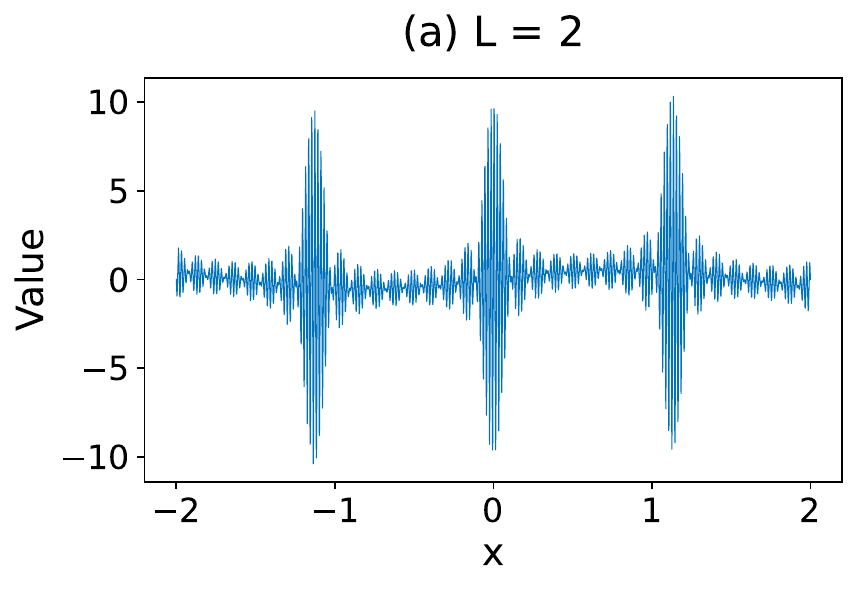}}  
\hspace{0.08\linewidth}
\subfloat{\includegraphics[width=0.35\linewidth]{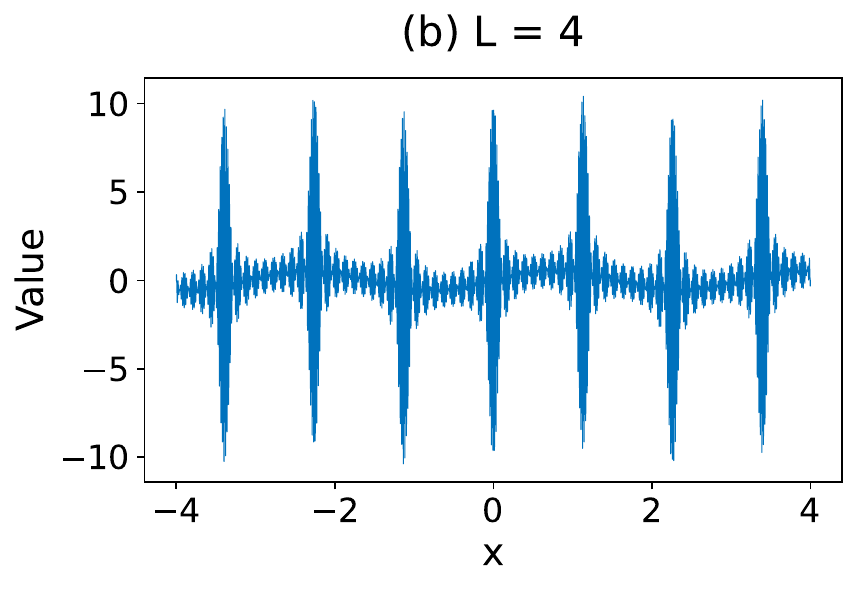}} \\
\subfloat{\includegraphics[width=0.35\linewidth]{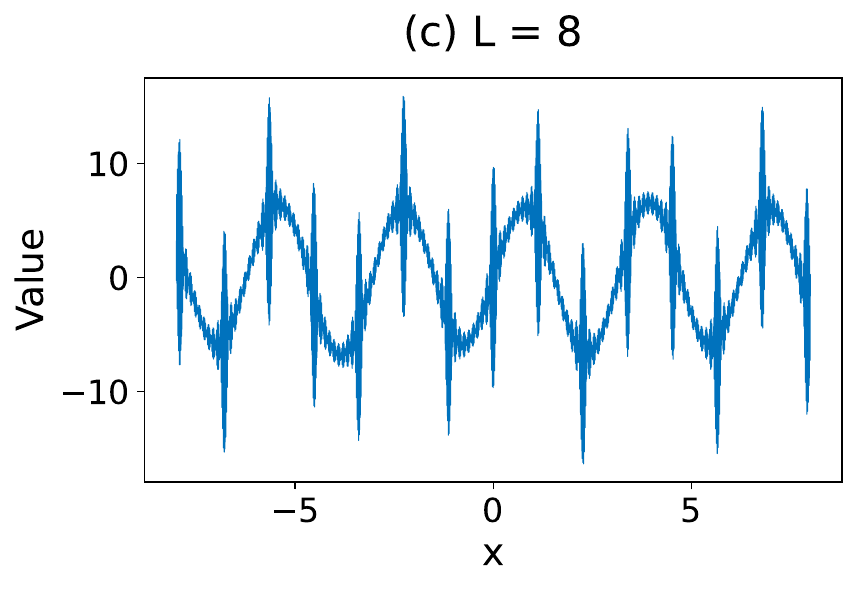}}  
\hspace{0.08\linewidth}
\subfloat{\includegraphics[width=0.35\linewidth]{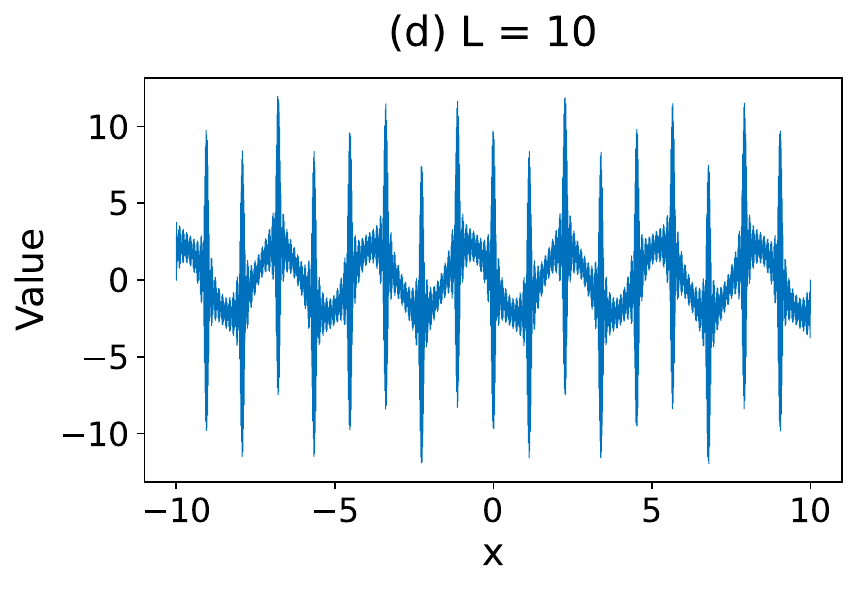}} \\
\caption{Characteristic solutions of the Helmholtz equation in spatial space for different domain lengths $L$} 
\label{helm_50_o}
\end{figure}

\begin{figure}[h]  
\centering  
\subfloat{\includegraphics[width=0.35\linewidth]{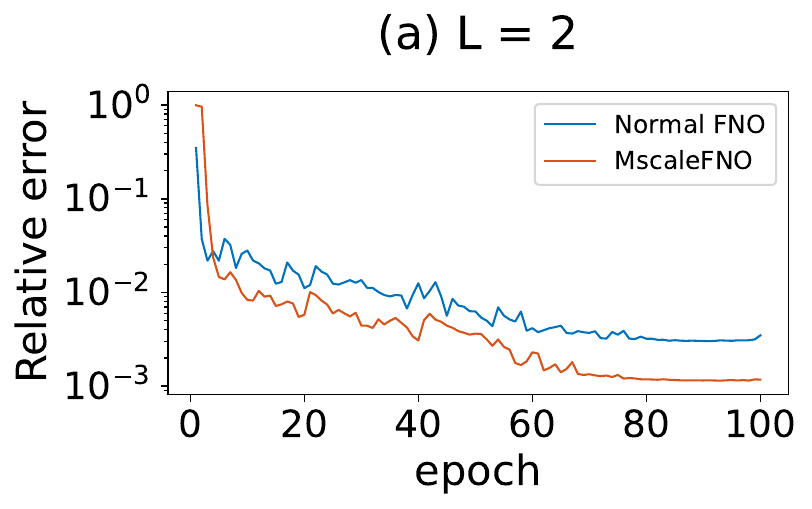}}  
\hspace{0.08\linewidth}
\subfloat{\includegraphics[width=0.35\linewidth]{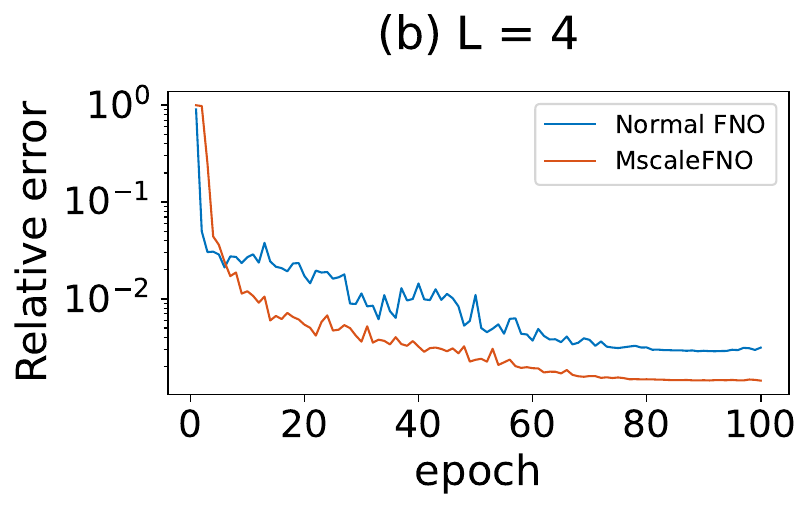}} \\
\subfloat{\includegraphics[width=0.35\linewidth]{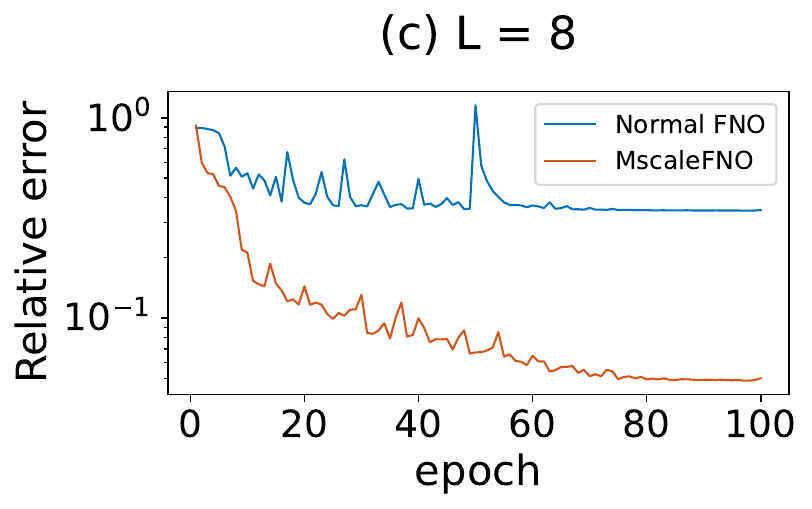}}  
\hspace{0.08\linewidth}
\subfloat{\includegraphics[width=0.35\linewidth]{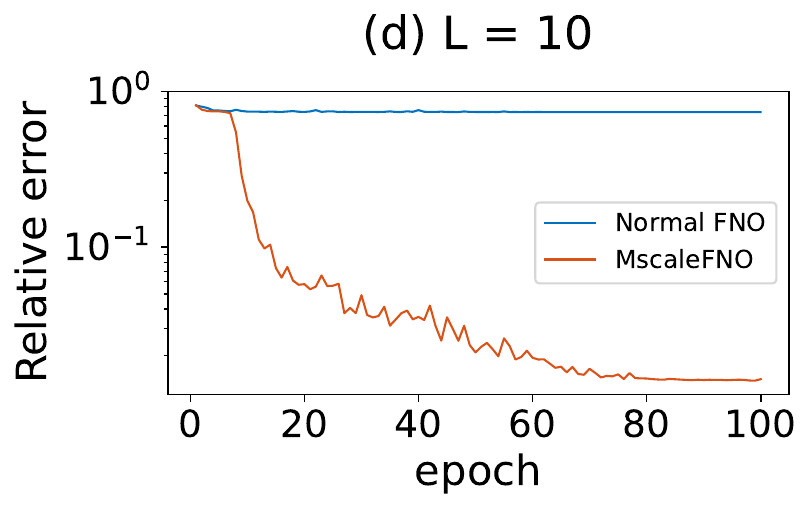}} \\
\caption{Error curves under different values of $L$ (Epoch=100)} 
\label{helm_50_ec}
\end{figure}

Fig.~\ref{helm_50_o} displays characteristic solutions of the Helmholtz equation for varying domain lengths $L \in \{2,4,8,10\}$, demonstrating the increasingly complex wave patterns that emerge as the spatial domain expands.
Fig.~\ref{helm_50_ec} presents the relative testing error curves for both normal FNO and MscaleFNO under different domain lengths $L \in \{2,4,8,10\}$ over 100 training epochs. As the length $L$ of the solving domain increases, corresponding to scattering problems in higher frequency regimes, the normal FNO (blue curves) exhibits progressively deteriorating performance, with the relative testing error reaching approximately 0.7 when $L=10$. In contrast, the MscaleFNO (orange curves) demonstrates robust performance across all domain lengths, not only consistently outperforming the normal FNO but also maintaining a relative error accuracy better than $O(10^{-2})$ throughout different spatial scales.

\begin{figure}[h]  
\centering  
\subfloat{\includegraphics[width=0.45\linewidth]{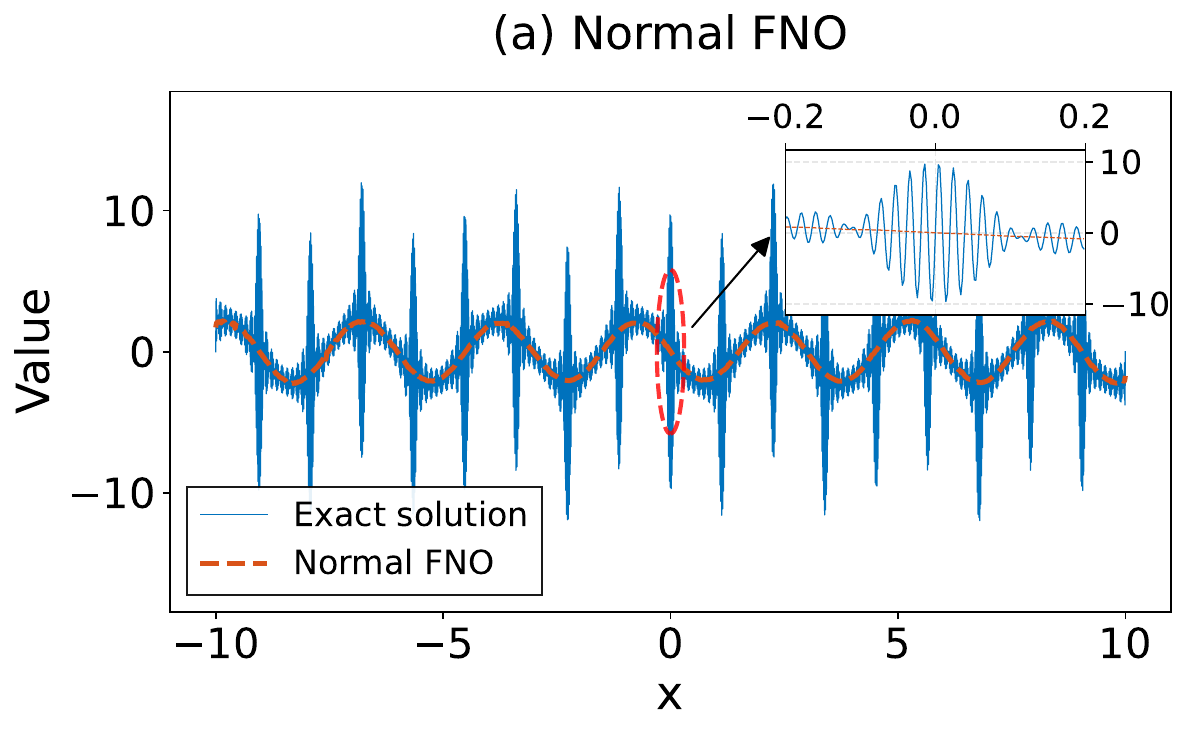}}  
\hspace{0.04\linewidth}
\subfloat{\includegraphics[width=0.45\linewidth]{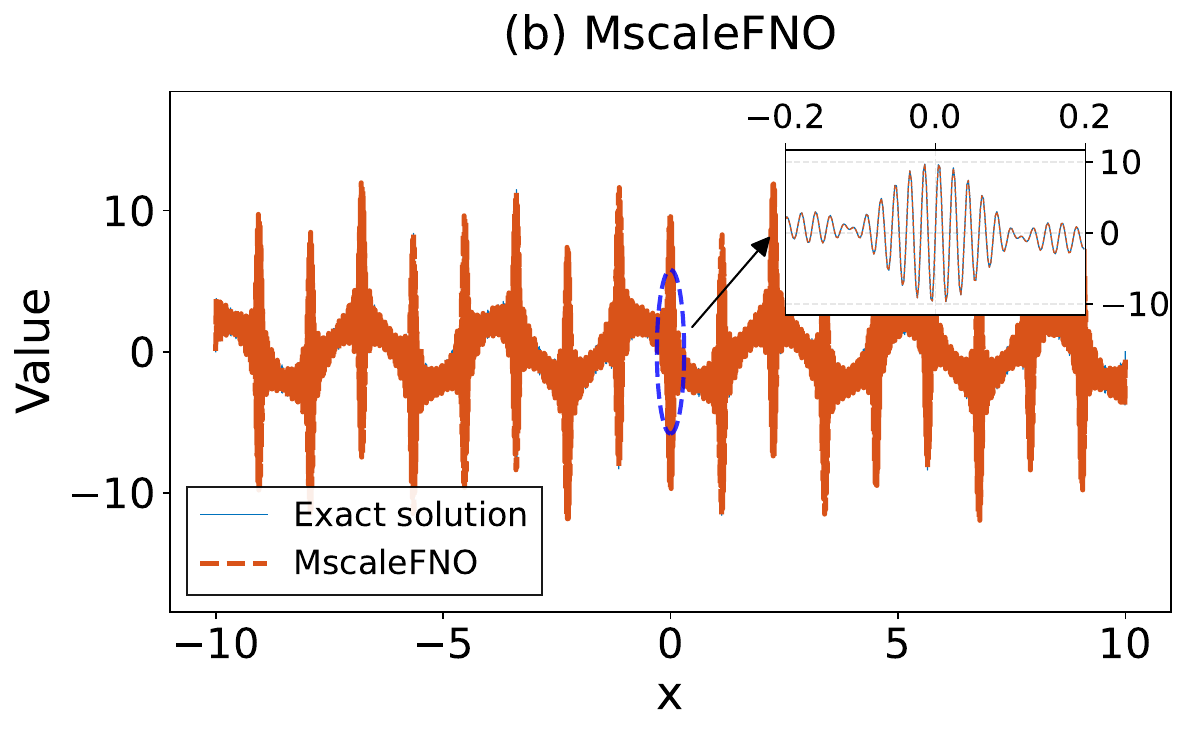}}  
\caption{ $L = 10$: Predicted solution by normal FNO (left) and MscaleFNO (right) with zoomed-in inset for $x\in[-0.2,0.2]$ } 
\label{helm_w50_L10_fitting_spatial}  
\end{figure}

\begin{figure}[h]  
\centering  
\subfloat{\includegraphics[width=0.45\linewidth]{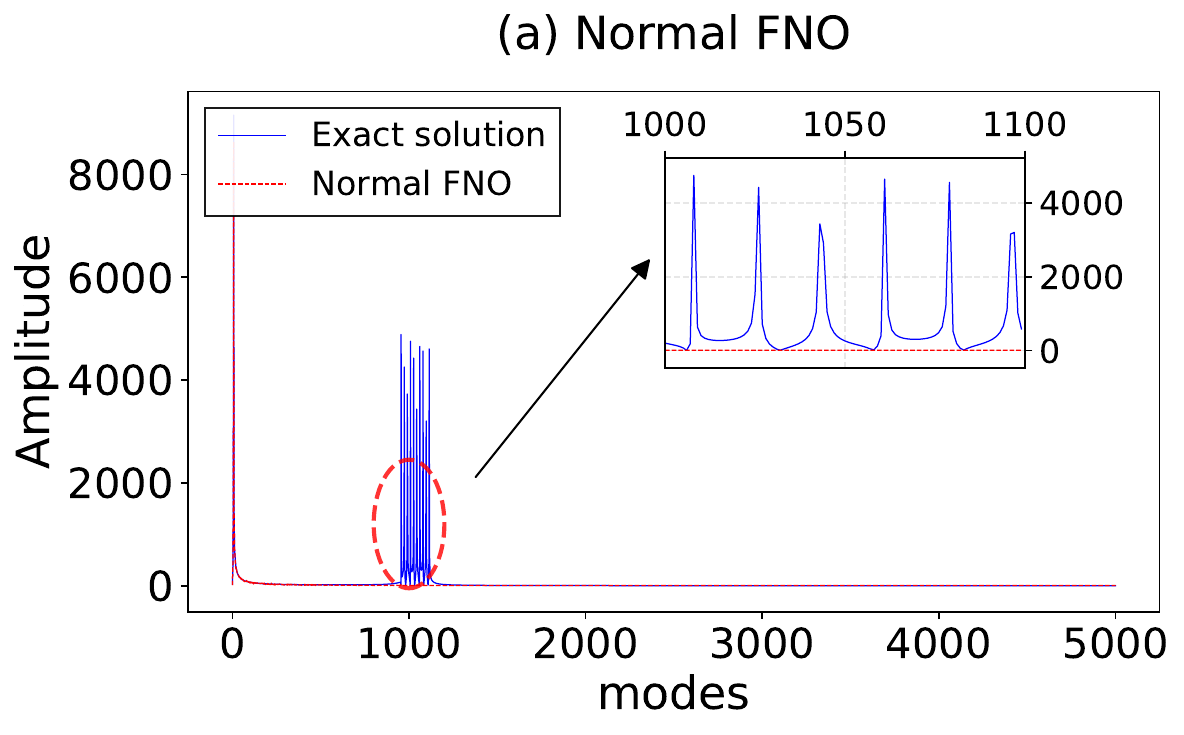}} 
\hspace{0.04\linewidth}
\subfloat{\includegraphics[width=0.45\linewidth]{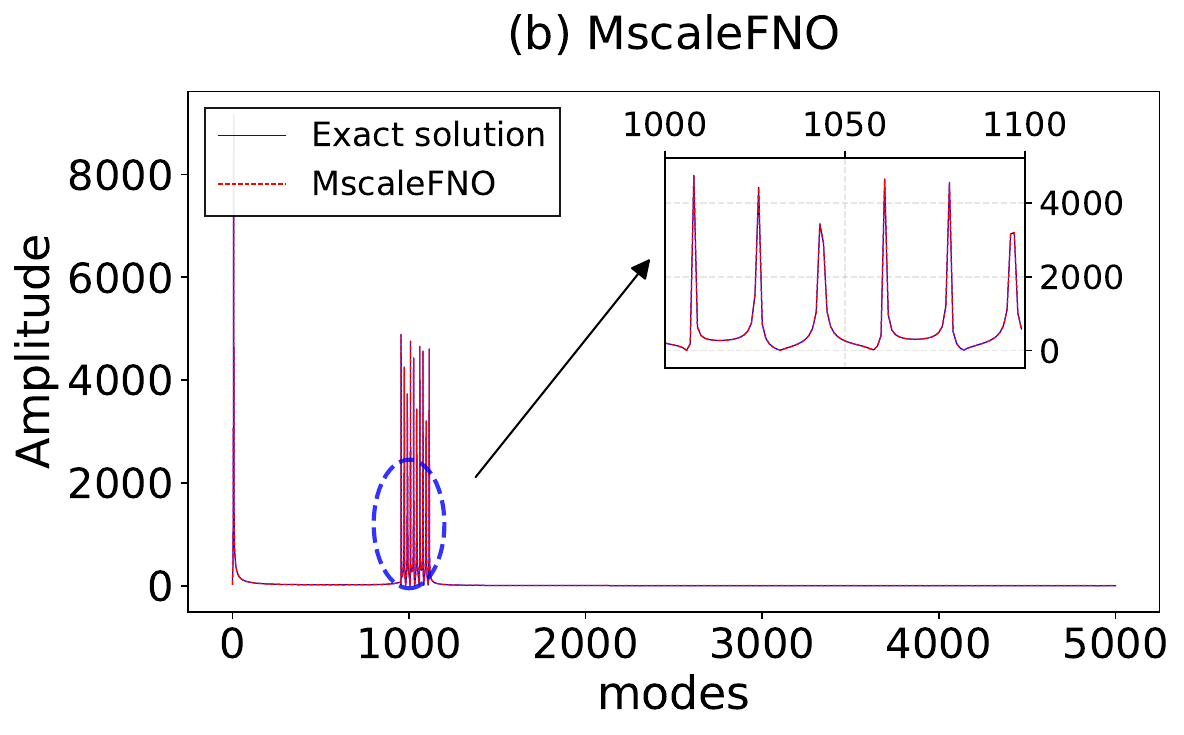}}  
\caption{$L = 10$: DFT of predicted solution by normal FNO (left) and MscaleFNO (right) with zoomed-in inset for modes $\in[1000,1100]$} 
\label{helm_w50_L10_fitting_spec}  
\end{figure}
Fig.~\ref{helm_w50_L10_fitting_spatial} and Fig.~\ref{helm_w50_L10_fitting_spec} present the prediction results for the case of $L=10$. In the physical domain (Fig.~\ref{helm_w50_L10_fitting_spatial}), normal FNO captures the general wave patterns but fails to get the high-frequency details, while MscaleFNO accurately reproduces the full solution structure including high-frequency oscillations.  
The Fourier analysis (Fig.~\ref{helm_w50_L10_fitting_spec}) supports these observations. Normal FNO preserves low-frequency components but shows significant distortion at high frequencies, whereas MscaleFNO accurately reconstructs the entire spectrum, as evidenced by the zoomed-in inset (frequencies 1000-1100).

\begin{figure}[h]  
\centering  
\subfloat{\includegraphics[width=0.45\linewidth]{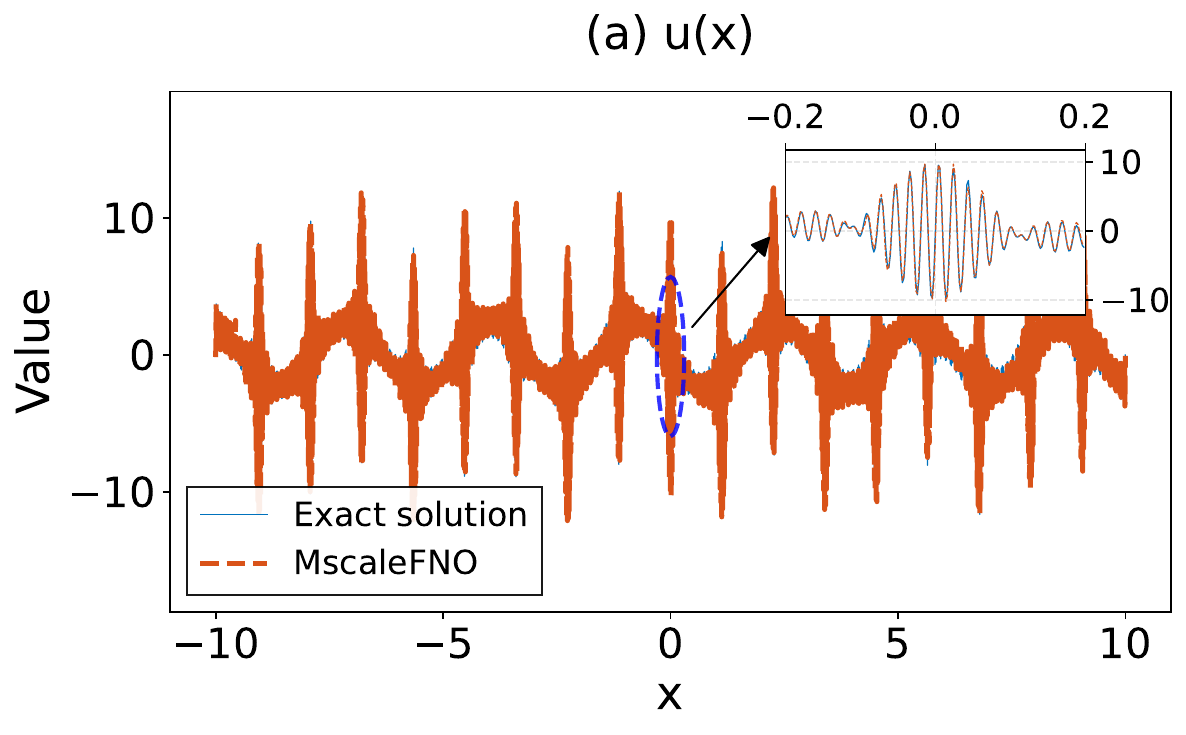}}  
\hspace{0.04\linewidth}
\subfloat{\includegraphics[width=0.45\linewidth]{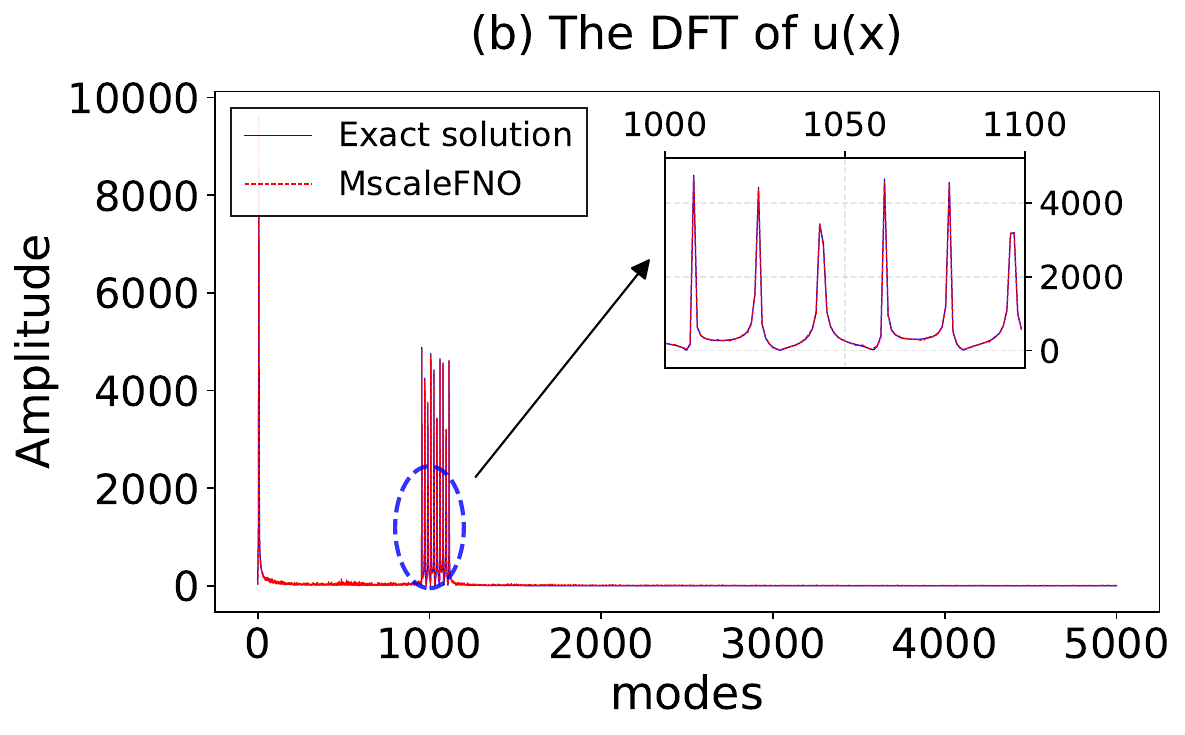}}  
\caption{$L = 10$: (a) Predicted solution of MscaleFNO against exact solution with zoomed-in inset for $x\in[-0.2,0.2]$ and (b) The DFT of $u(x)$ with zoomed-in inset for modes $\in[1000,1100]$ } 
\label{helm_w50_test_fitting}  
\end{figure}

\begin{eg}
To evaluate the generalization capability of the MscaleFNO, we perform numerical experiments at $L=10$ on test samples drawn from distributions distinct from the training data,
$$\eta(x)=\sum_{n=1}^{50}a_n \sin(k_n x^3)+b_n \cos(l_n x^2),\; \omega(x) = \frac{\eta(x)}{\max_x\big\{\eta(x)\big\}},$$  
where $k_n\sim \hbox{rand}(0,30),l_n\sim \hbox{rand}(40,60)$ and $f$ is the same as Eq.~\eqref{e11}.  
The equation solving and data acquisition procedures remain consistent with the previous setup. 
Fig~\ref{helm_w50_test_fitting_normal} shows that normal FNO fails on unseen test functions, with significant prediction errors in both solution amplitude and oscillation patterns. In contrast, MscaleFNO Fig~\ref{helm_w50_test_fitting}demonstrates robust prediction capability. Even with a change of test function form not seen during the training, the model accurately captures high-frequency oscillations over intervel $[-10,10]$ and precisely predicts both low and high-frequency modes in the Fourier spectrum (modes 1000-1100). This robust spectral reconstruction ability, combined with excellent generalization beyond training distribution, suggests MscaleFNO's potential for accurately approximating arbitrary continuous functions.
\begin{figure}[h]  
\centering  
\subfloat{\includegraphics[width=0.45\linewidth]{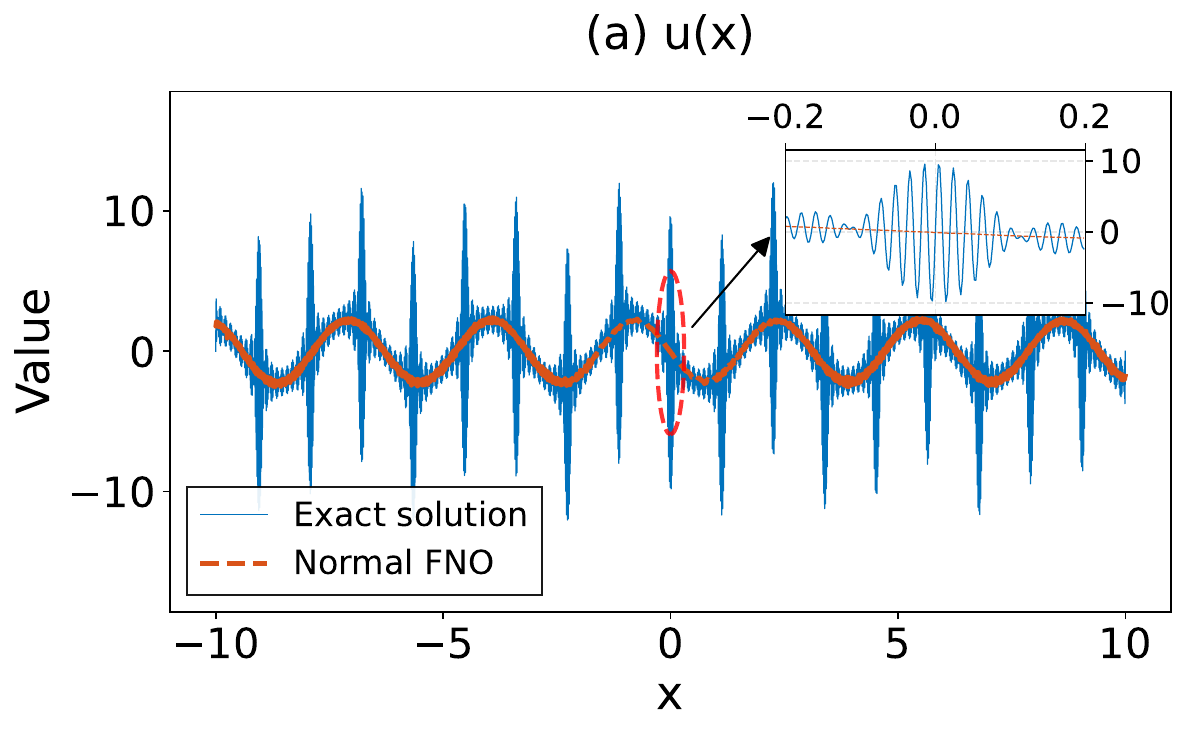}}  
\hspace{0.04\linewidth}
\subfloat{\includegraphics[width=0.45\linewidth]{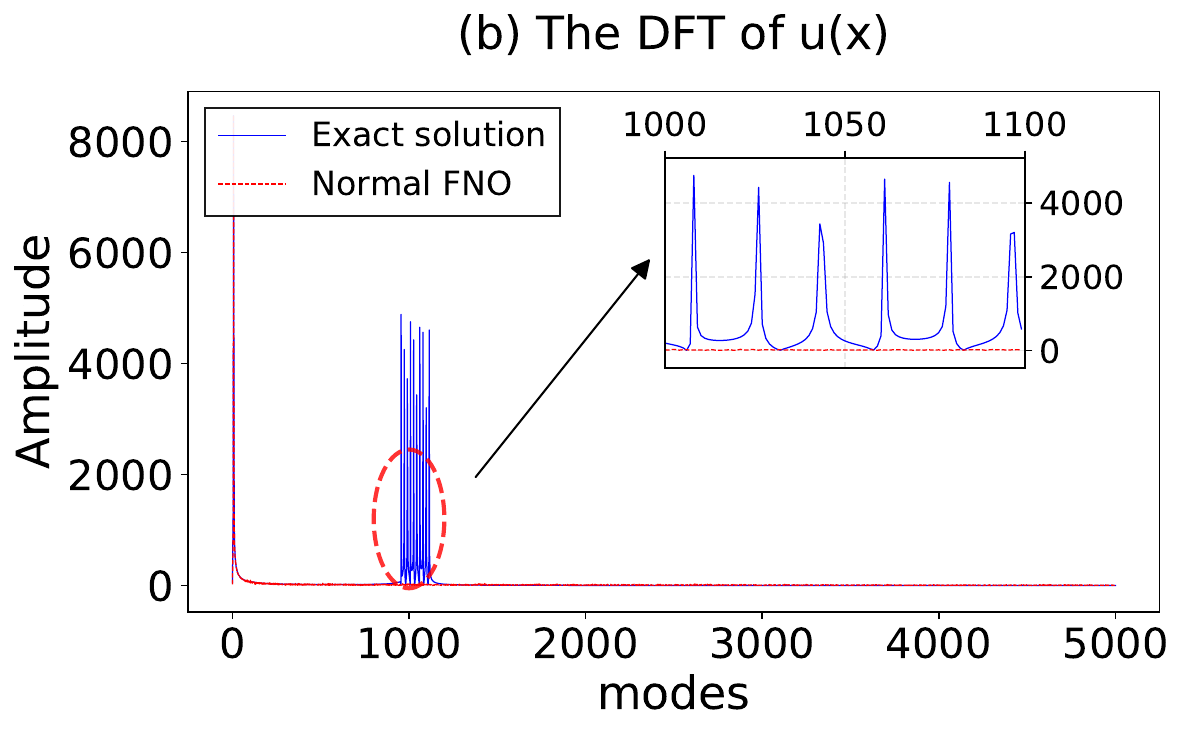}}  
\caption{$L = 10$: (a) Predicted solution of normal FNO against exact solution with zoomed-in inset for $x\in[-0.2,0.2]$ and (b) The DFT of $u(x)$ with zoomed-in inset for modes $\in[1000,1100]$ } 
\label{helm_w50_test_fitting_normal}  
\end{figure}

\end{eg}  

\section{Conclusion}
We have proposed a multi-scale Fourier neural operator to improve the spectral representation of the neural network in learning mapping between oscillatory functions, such as the one between material properties of a scatterer and its scattering field, especially in high-frequency regime. Numerical results have shown significant improvements of the MscaleFNO over the normal FNO. Future work will include applying the MscaleFNO for higher dimensional Helmhotlz equations as well as to solve inverse medium problems in high frequency wave scattering.

\section*{Acknowledgement}
Z. You and Z. Xu are financially supported by the National Natural Science Foundation of China (Grant Nos. 12325113 and 12426304). W. Cai is supported by the Clements Chair of Applied Mathematics at SMU.


\end{document}